\documentclass[a4paper,11pt]{article}
\usepackage{amssymb}

\newcommand{\R}{\mathbb{R}}

\newcommand{\cuad}{{\sqcap\kern-.68em\sqcup}}

\newcommand{\norm}[1]{\|#1\|}

\newtheorem{theorem}{Theorem}[section]
\newtheorem{proposition}{Proposition}[section]

\newtheorem{lemma}{Lemma}[section]

\newtheorem{remark}{Remark}[section]
\newcommand{\bremark}{\begin{remark} \em}
\newcommand{\eremark}{\end{remark} }

\begin{document}

\begin{center}{\bf  \large  On nonhomogeneous  elliptic  equations  with \\[2mm]

  the Hardy-Leray potentials  }\medskip
%%%%%%%%%%%%%%%%%%%%%%%%%%%%%%%%%%%%%%%%%%%%%%%%%%%%%%%%%%%%%%%%%%%%%%
%%%%%%%%%%%%%%%%%%%%%%%%%%%%%%%%%%%%%%%%%%%%%%%%%%%%%%%%%%%%%%%%%%%%%%

 \bigskip\medskip

 {\small  Huyuan Chen\footnote{chenhuyuan@yeah.net}, \quad Alexander Quaas\footnote{alexander.quaas@usm.cl}\quad and\quad  Feng Zhou\footnote{ fzhou@math.ecnu.edu.cn}
}
 \bigskip\medskip

{\small  $ ^1$Department of Mathematics, Jiangxi Normal University,\\
Nanchang, Jiangxi 330022, PR China\\[3mm]

$ ^2$Departamento de Matem\'{a}tica, Universidad T\'{e}cnica Federico,
Santa Mar\'{i}a\\ Casilla V-110, Avda. Espa\~{n}a 1680,
Valpara\'{i}so, Chile \\[3mm]

$ ^3$Center for PDEs and Department of Mathematics, East China Normal University,\\
 Shanghai, 200241, PR China
} \\[6mm]

\begin{abstract}
In this paper,  we present some suitable distributional identities  of  solutions for nonhomogeneous elliptic equations involving the Hardy-Leray potentials and
study qualitative properties of the solutions to the corresponding nonhomogeneous problems  in this  distributional sense. We address some applications on the nonexistence of
some nonhomogeneous problems with the Hardy-Leray potentials and  the nonexistence principle eigenvalue with some indefinite potentials.

\end{abstract}

\end{center}
%\tableofcontents \vspace{1mm}
  \noindent {\small {\bf Keywords}:   Distributional identity,  Hardy-Leray potential, Isolated singularity.}\vspace{1mm}

\noindent {\small {\bf MSC2010}:     35B40,     35J99. }

\vspace{2mm}
%%%%%%%%%%%%%%%%%%%%%%%%%%%%%%%%%%%%%%%%%%%%%%%%%%%%%%%%%%%%%%%%%%%%%%%%%%%%%%%%%%%%%%%%%%%%%%%%%%%%%%%%%%%%%%%%%%%%%%%%%%
%%%%%%%%%%%%%%%%%%%%%%%%%%%%%%%%%%%%%%%%%%%%%%%%%%%%%%%%%%%%%%%%%%%%%%%%%%%%%%%%%%%%%%%%%%%%%%%%%%%%%%%%%%%%%%%%%%%%%%%%%%

\setcounter{equation}{0}
\section{Introduction}

Throughout this paper, we assume that $\Omega$ is an open domain containing the origin in $\R^N$ with $N\ge 2$, $\mu\ge \mu_0:=-\frac{(N-2)^2}{4}$   and $\mathcal{L}_\mu:= -\Delta   +\frac{\mu}{|x|^{2 }}.$  The elliptic operator  $\mathcal{L}_\mu$
involves the inverse square potential, which is also called as Hardy potential,   and related semilinear elliptic problems have been done mainly by  variational methods \cite{CL,DD,D,FF},
 due to the Hardy type inequalities, see the references \cite{ACR,BM,BV1,GP}.
A direct extension of definition for the distributional solution to the Hardy problem $\mathcal{L}_\mu u=g\ {\rm in}\  \Omega$
may be proposed as
\begin{equation}\label{d1}
 \int_\Omega u\mathcal{L}_\mu \xi\, dx=\int_\Omega g\xi\, dx,\quad\forall\,\xi\in C^\infty_c(\Omega),
\end{equation}
where  $g$ is a nonlinearity of $x$ and $u$. When $N\ge3$ and $\mu\in[\mu_0,0)$, \cite{BP,D} make use of this notation of distributional solution to
show the existence of distributional solutions of (\ref{d1}) for particular nonlinearity depending on $u$ under some restriction for $\mu$.
Later on, the authors in \cite{FM1} considered the related Hardy problem in the   distributional sense of (\ref{d1}) replaced $C^\infty_c(\Omega)$
by $C^\infty_c(\Omega\setminus\{0\})$ in some range of $\mu$.
 Chaudhuri and C\^irstea in \cite{CC}, C\^irstea in \cite{C}  classified the isolated singular classical solution of $\mathcal{L}_\lambda u+b(x)h(u)=0$ in $\Omega\setminus \{0\}$, where
both $b$ and $h$ consist of regularly varying and slowly varying parts (see their definitions in \S1.2.2 of \cite{C}). There a solution is considered as a $C^1(\Omega \setminus \{0\})$-solution in the sense of distributions in $\Omega \setminus \{0\}$, that is,
\begin{equation}\label{d11}
\int_\Omega \nabla u \nabla \varphi dx - \int_\Omega \frac{\lambda}{|x|^2} u \varphi dx + \int_\Omega b(x)h(u) \varphi dx =0, \quad\forall\,\varphi \in C^1_c(\Omega\setminus \{0\})
\end{equation}
holds. Recall that  a solution $u$ has   a "removable" singularity at the origin  if $u$ can be extended to a $C^1$-solution in $\mathcal D'(\Omega)$. Another important
attempt is done in \cite{D} to  consider the consider the solutions of semilinear Hardy equation in classical sense, in variational sense and  in a distributional sense.
 This subject has been extended into the elliptic equation with Hardy
potential with singularities on  the  boundary,  for instance, \cite{BMM,FM2,GT,MM,MN,KV} and the references therein.

However, we observe that the distributional identity (\ref{d11}) loses the information at the origin in the case that the test functions  are in $C^\infty_c(\Omega\setminus\{0\})$,
and on the other hand, it does not work  even for the fundamental solutions by the test functions in $C^\infty_c(\Omega)$ if $\mu>0$ large enough when we consider the distributional identity (\ref{d1}), because of high singularity at the origin of the corresponding fundamental solutions.     Furthermore,
it loses the huge convenience in dealing with semi-linear elliptic equations involving the Laplacian  from the expression of singularities by Dirac mass, for instances \cite{BL,L,NS}, where
the authors obtained classical solutions of $-\Delta u=f$ with precise isolated singularities  by dealing with $-\Delta u=f+k\delta_0$ in models: $f$ being nonhomogeneous term,  $f(u)=u^p$ or $f(u)=u^p-u$.

Therefore, it is important for theoretical advances and also for application's interests to obtain an improved version of distributional idnetity (\ref{d1}).  To this end, we start our analysis form the fundamental solutions of Hardy operators.   It is known that when $\mu\ge\mu_0$,  the problem
\begin{equation}\label{eq 1.1}
\mathcal{L}_\mu u= 0\quad {\rm in}\quad  \R^N\setminus \{0\}
\end{equation}
has two branches of radially symmetric solutions with the explicit formulas that
\begin{equation}\label{1.1}
 \Phi_\mu(x)=\left\{\arraycolsep=1pt
\begin{array}{lll}
 |x|^{\tau_-(\mu)}\quad
   &{\rm if}\quad \mu>\mu_0,\\[1mm]
 \phantom{   }
-|x|^{\tau_-(\mu)}\ln|x| \quad &{\rm   if}\quad \mu=\mu_0
 \end{array}
 \right.\quad \quad{\rm and}\quad \Gamma_\mu(x)=|x|^{ \tau_+(\mu)},
\end{equation}
where
$$
 \tau_-(\mu)=-\frac{N-2}2-\sqrt{\mu-\mu_0}\;\;\;{\rm and}\;\;\;  \tau_+(\mu)=-\frac{ N-2}2+\sqrt{\mu-\mu_0}.
$$
Here $\tau_-(\mu)$ and $\tau_+(\mu)$ are two roots of $\mu-\tau(\tau+N-2)=0$. Moreover, $ \Phi_\mu$ is a regular solution of
(\ref{eq 1.1}) in the sense that $\mu |\cdot|^{-2} \Phi_{\mu}(\cdot) \in L^1_{loc}(\R^N)$ and we have that
\begin{equation}\label{eq 1.100}
\mathcal{L}_\mu \Phi_{\mu}= 0\quad {\rm in}\quad  \mathcal D'(\R^N).
\end{equation}
Remark that the  mentioned fundamental solution of (\ref{eq 1.1}) with $\mu=0$ is
$$\Phi_0(x)= \left\{\arraycolsep=1pt
\begin{array}{lll}
 |x|^{2-N}\quad
   &{\rm if}\quad N\ge  3,\\[1mm]
 \phantom{   }
- \ln|x| \quad  &{\rm   if}\quad N=2
 \end{array}
 \right.$$
and      $\Phi_0$ is a distributional solution of
\begin{equation}\label{eq003}
\mathcal{L}_0  u =c_{N}\delta_0 \quad   \rm{in}\quad\R^N,
\end{equation}
where
$$
 c_N=\left\{\arraycolsep=1pt
\begin{array}{lll}
 (N-2)|\mathcal{S}^{N-1}|\quad
   &{\rm if}\quad N\ge 3,\\[1.5mm]
 \phantom{   }
2\pi \quad  &{\rm  if}\quad N=2
 \end{array}
 \right.
$$
and $\mathcal{S}^{N-1}$ is the sphere of  the unit ball in $\R^N$ and $|\mathcal{S}^{N-1}|$ is the area of the unit sphere. We observe that the second branch of
fundamental solution $\Gamma_0\equiv1$ for $\mu=0$, which is always omitted in the study of semilinear elliptic equation. Without special stating, we always assume that $\mu\ge \mu_0$.

 Our motivation of this paper is to provide   suitable distributional identities for Hardy problems and to  answer the basic question whether it is able to express the isolated singularities by the Dirac mass, in particular,  to find out some  distributional identity to distinguish the fundamental solutions of $\mathcal{L}_\mu$.
 %In this article, we find out a creative way to express the related singularities by the Dirac mass.

Our first result on the distributional identities of the fundamental solutions can be stated as follows:
 \begin{theorem}\label{teo 1}
Let $d\mu(x) =\Gamma_\mu(x) dx$ and
\begin{equation}\label{L}
 \mathcal{L}^*_\mu=-\Delta -\frac{2 \tau_+(\mu) }{|x|^2}\,x\cdot\nabla.
\end{equation}
Then
\begin{equation}\label{1.2}
 \int_{\R^N}\Phi_\mu   \mathcal{L}^*_\mu(\xi) \, d\mu  =c_\mu\xi(0),\quad \forall\, \xi\in C^{1.1}_c(\R^N),
\end{equation}
where
 \begin{equation}\label{cmu}
 c_\mu=\left\{\arraycolsep=1pt
\begin{array}{lll}
 2\sqrt{\mu-\mu_0}\,|\mathcal{S}^{N-1}|\quad
   &{\rm if}\quad \mu>\mu_0,\\[1.5mm]
 \phantom{   }
|\mathcal{S}^{N-1}| \quad  &{\rm  if}\quad \mu=\mu_0.
 \end{array}
 \right.
 \end{equation}

\end{theorem}

Here and in what follows, we always take the notation
$$d\mu(x)=\Gamma_\mu(x) dx\quad{(\rm\, or}\ d\mu=\Gamma_\mu dx).$$

The distributional identity (\ref{1.2}) is derived from the observation that $\mathcal{L}_\mu(\Gamma_\mu\xi)=\Gamma_\mu \mathcal{L}_\mu^*(\xi)$.
In the $d\mu$-distributional sense, $\mathcal{L}_\mu^*$ is the duality operator of $\mathcal{L}_\mu$   and we may  say that the fundamental solution $\Phi_\mu$ is a   $d\mu-$distributional solution of
$$%\begin{equation}\label{eq001}
 \mathcal{L}_\mu \Phi_\mu  =c_\mu\delta_0 \quad   \rm{in}\quad\mathcal D'(\R^N),
$$%\end{equation}
which coincides (\ref{eq003}) when $\mu=0$, since $d\mu(x)=dx$ in this case and $\mathcal{L}_0^*=\mathcal{L}_0=-\Delta$. Furthermore, an alternative form of $\mathcal{L}_\mu^*$
could be given as
 $$\mathcal{L}_\mu^*(\xi)=- \Gamma_\mu^{-1}{\rm div}(\Gamma_\mu \nabla \xi),$$
 which plays an essential role in the study of  heat kernel for the Hardy operator, see \cite{FMT,MT}. We remark that
 when $\mu=0$, $c_\mu$ coincides with the coefficient $c_N$, while for $\mu=\mu_0$  there is a jump for the parameter $c_\mu$.

It is known  that the operator $\mathcal{L}_\mu$ with   $\mu\in[\mu_0,+\infty)\setminus\{0\}$, is self-adjoint in the normal distributional sense. However,  identity (\ref{1.2}) gives a new distributional sense with respect to a specific weighted measure,
which  makes the operator $\mathcal{L}_\mu$ not self-adjoint.
 To the best of our knowledge, identity (\ref{1.2}) is new.   \smallskip

% For the distributional identities (\ref{1.2}) with $\mu=0$ see \cite{V} and the references therein.  We remark that  the identity (\ref{1.2}) is equivalent to
%$$\int_{\R^N}\Phi_\mu \mathcal{L}_\mu(\Gamma_\mu\xi) dx :=\lim_{r\to0^+} \int_{\R^N\setminus B_r(0)}\Phi_\mu \mathcal{L}_\mu(\Gamma_\mu\xi) dx=c_\mu\xi(0).$$
% We note that from Theorem \ref{teo 1} it indicates that two branches of the fundamental solutions cooperate and generate the Dirac mass in the distributional sense.

  We next  continue to extend this identity for bounded domains and apply it  for the classification of
 isolated singularities of nonhomogeneous Hardy problem.   Given a bounded $C^2$ domain $\Omega$ containing the origin, the Hardy problem
 \begin{equation}\label{eq 1.1b}
 \arraycolsep=1pt\left\{
\begin{array}{lll}
 \displaystyle  \mathcal{L}_\mu u= 0\qquad
   {\rm in}\quad  {\Omega}\setminus \{0\},\\[1mm]
 \phantom{  L_\mu \, }
 \displaystyle  u= 0\qquad  {\rm   on}\quad \partial{\Omega},\\[1mm]
 \phantom{   }
  \displaystyle \lim_{x\to0}u(x)\Phi_\mu^{-1}(x)=1
 \end{array}\right.
\end{equation}
 has  solution $G_\mu$ with isolated singularity as $\Phi_\mu$ at origin, i.e.
\begin{equation}\label{1.3}
 \lim_{|x|\to0} \frac{G_\mu(x)}{\Phi_\mu(x)}=1.
\end{equation}
%where, throughout this paper, we assume that $\Omega$ is a bounded domain in $\R^N$ containing the origin.
Then we have the following distributional identity.

\begin{theorem}\label{teo 2}
Let   $G_\mu$ be the solution of (\ref{eq 1.1b}) verifying (\ref{1.3}),  then
\begin{equation}\label{1.2b}
 \int_{\Omega}G_\mu  \mathcal{L}_\mu^*(\xi)\, d\mu  =c_\mu \xi(0),\quad \forall\,\xi\in  C^{1.1}_0(\Omega).
\end{equation}
\end{theorem}

%We emphasis that the distributional identity (\ref{1.2b}) holds for the measure $d\mu = \Gamma_\mu dx$,  then the function $G_\mu$ could be considered as a

A deeper knowledge of   distributional identities   allows us to draw a   complete picture of the existence, non-existence and the singularities
for  the  nonhomogeneous problem
 \begin{equation}\label{eq 1.1f}
 \arraycolsep=1pt\left\{
\begin{array}{lll}
 \displaystyle   \mathcal{L}_\mu u= f\qquad
   {\rm in}\quad  {\Omega}\setminus \{0\},\\[1.5mm]
 \phantom{   L_\mu   }
 \displaystyle  u= 0\qquad  {\rm   on}\quad \partial{\Omega}£¬
 \end{array}\right.
\end{equation}
where $f:\bar\Omega\setminus\{0\}\mapsto \R$ is a H\"{o}lder continuous locally in $\bar\Omega\setminus\{0\}$.

 Motivated by (\ref{1.2}) and (\ref{1.2b}),  we shall classify the isolated singularities  of (\ref{eq 1.1f})  by building the connection with
\begin{equation}\label{eq 1.2}
 \arraycolsep=1pt\left\{
\begin{array}{lll}
 \displaystyle   \mathcal{L}_\mu u= f+c_\mu k\delta_0\qquad
   &{\rm in}\quad  {\Omega},\\[1.5mm]
 \phantom{   L_\mu   }
 \displaystyle  u= 0\qquad  &{\rm   on}\quad \partial{\Omega},
 \end{array}\right.
\end{equation}
in the $d\mu$-distributional sense that $u\in L^1(\Omega, d\mu)$,
 \begin{equation}\label{1.2f}
 \int_{\Omega}u  \mathcal{L}_\mu^*(\xi)\, d\mu  = \int_{\Omega} f  \xi\, d\mu +c_\mu k\xi(0),\quad\forall\, \xi\in   C^{1.1}_0(\Omega),
\end{equation}
where $k\in\R$.
Note that the solution  $G_\mu$ of   (\ref{eq 1.1b}) verifying (\ref{1.3}),  from Theorem \ref{teo 2}, is a $d\mu$-distributional solution of
 $$  %\label{eq 1.1bd}
 \arraycolsep=1pt\left\{
 \begin{array}{lll}
  \displaystyle   \mathcal{L}_\mu u = c_\mu \delta_0\qquad
   &{\rm in}\quad  {\Omega},\\[1mm]
 \phantom{  L_\mu   }
 \displaystyle  u= 0\qquad  &{\rm   on}\quad \partial{\Omega}.
 \end{array}\right.
 $$

For nonhomogeneous problem (\ref{eq 1.1f}), we have the following results.

\begin{theorem}\label{teo 3}
Let    $f$ be a function in $C^\gamma_{loc}(\overline{\Omega}\setminus \{0\})$ for some $\gamma\in(0,1)$.

$(i)$ Assume that
\begin{equation}\label{f1}
f \in L^1(\Omega, d\mu)\qquad {\rm i.\,e.}\quad  \int_{\Omega} |f|\,   d\mu <+\infty,
\end{equation}
then  for any $k\in\R$, problem (\ref{eq 1.2}) admits a unique weak solution  $u_k$, which is a classical solution of   problem (\ref{eq 1.1f}).
Furthermore, if assume more that
\begin{equation}\label{4.1-3}
 \lim_{|x|\to0}f(x)|x|^{2-\tau_-(\mu)}=0,
\end{equation}
then  we have the asymptotic behavior
\begin{equation}\label{singular}
\lim_{x\to0}u_k(x)\Phi_\mu^{-1}(x)=k.
\end{equation}

$(ii)$ Assume that $f$ verifies (\ref{f1})  and $u$ is a nonnegative solution of (\ref{eq 1.1f}), then  $u$ is a weak solution of  (\ref{eq 1.2}) for some $k\ge0$.\smallskip

$(iii)$ Assume that $f\ge0$ and
\begin{equation}\label{f2}
 \lim_{r\to0^+} \int_{\Omega\setminus B_r(0)} f\, d \mu  =+\infty,
\end{equation}
then problem (\ref{eq 1.1f}) has no nonnegative solutions.
\end{theorem}

Note that Theorem \ref{teo 3} part $(i)$ and $(ii)$ build a one to one connection between problems (\ref{eq 1.1f}) and (\ref{eq 1.2}) under the assumption  (\ref{f1}). Together with Theorem \ref{teo 3} part $(iii)$,  the assumptions of (\ref{f1}) is sharp  for the existence of weak solution to (\ref{eq 1.1f}).
Assumption (\ref{f1}) could also be seen in  \cite{D}, where the authors studied the existence of isolated singular solutions of (\ref{eq 1.1f}) only when $N\ge 3$ and $\mu_0\le \mu<0$. Similar assumptions seem be initially proposed  in \cite{BG,BG1} to  deal with semilinear parabolic equation with Hardy potentials.

 It worth noting that our weak solution $u_k$ of (\ref{eq 1.2})  could be decomposed into $v_f+kG_\mu$,
where  $v_f$ is a unique solution of (\ref{1.2f}) with $k=0$ and $G_\mu$ is the solution of (\ref{eq 1.1b}) subject to (\ref{1.3}),  and the classification of the isolated singularities
of (\ref{eq 1.1f}) is performed by Schwartz Theorem in \cite{S}.

 Furthermore,  we enlarge the scope of the inhomogeneity into  $f\in L^1(\Omega,\, \rho d\mu)$, where $\rho(x)=dist(x,\partial\Omega)$ in Section 5.2 and  show the existence of   $d\mu$-distributional solution of (\ref{eq 1.2}) by approximating the inhomogeneity by regular functions.

By applying the nonexistence result of Theorem \ref{teo 3} $(iii)$, we obtain an Liouville theorem for the nonhomogeneous problem (\ref{eq 1.1f}) even for $\mu<\mu_0$  in Section 5.2.
Finally, our results could be applied to show the nonexistence of a principle eigenvalue problem with an indefinite potential with singularities as Hardy-Leray potentials.

The rest of the paper is organized as follows. In section 2, we  show   qualitative properties of the solutions to nonhomogeneous problem  with regular nonlinearity
and "removable" singularity $\lim_{x\to0}u(x)\Phi_\mu^{-1}(x)=0$. For the uniqueness, the basic tools are comparison principle in the classical sense and the basic tools Kato's inequality in the $d\mu$ weak sense.   Section 3 is devoted to build the distributional identity for the fundamental solutions in $\R^N$, and in bounded smooth domain
and we build the approximation for the fundamental solutions in the distributional sense of (\ref{1.2b}). Section 4 is addressed to
study the qualitative properties of the solutions of (\ref{eq 1.1f}) in the distributional sense  and  to prove Theorem \ref{teo 3}. Finally, we give a generalization of $d\mu$-distributional solution and an application to the nonexistence of the principle eigenvalue of Laplacian with zero Dirchlet boundary condition.

\setcounter{equation}{0}
\section{Preliminary}

\subsection{ Comparison Principle for $\mathcal L_{\mu}$ and weak singularity}
In this subsection, we introduce the Comparison Principle for the operator $\mathcal L_{\mu}$ which is one basic tool in our analysis. %as it is classical for the Laplacian operator.

\begin{lemma}\label{lm cp}%\cite{Ch}
Let $O$ be a bounded open set in $\R^N$, $L: O\times [0,+\infty)\to[0,+\infty)$ be a continuous function such that for any $x\in  O$,
$$L(x,s_1)\ge L(x,s_2)\quad {\rm if}\quad s_1\ge s_2,$$ then $\mathcal{L}_\mu+L$ with $\mu\ge \mu_0$ verifies the Comparison Principle,
that is, if
$$u,\,v\in C^{1,1}(O)\cap C(\bar O)$$ satisfying
$$\mathcal{L}_\mu u+ L(x,u)\ge \mathcal{L}_\mu v+ L(x,v) \quad {\rm in}\quad  O
\qquad{\rm and}\qquad   u\ge  v\quad {\rm on}\quad \partial O,$$
then
$$u\ge v\quad{\rm in}\quad  O.$$

\end{lemma}
{\bf Proof.} Let $w=u-v$ and $w_-=\min\{w,0\}$, then $w\ge 0$ on  $\partial O$ by the assumption that  $u\ge v$ on $\partial O$ and then
$w_-=0$ on $\partial O$.
We will prove that $w_-\equiv 0$.
If $ O_-:=\{x\in O:\, w(x)<0\}$ is not empty, then it is a bounded
$C^{1,1}$ domain in $ O$.  From Hardy inequality with $\mu\ge \mu_0$, for some $C>0$, there holds,
\begin{eqnarray*}
 0  &=& \int_{ O_-}(-\Delta  w_- +\frac{\mu}{|x|^{2 }}  w_- ) w_- dx+\int_{ O_-}[L(x,u)-L(x,v)](u-v)_-\,dx \\ &\ge&  \int_{ O_-} \left(|\nabla w_-|^2 +\frac{\mu}{|x|^{2 }}  w_-^2\right) dx
 \ge  C \int_{ O_-}w_-^2 dx,
\end{eqnarray*}
then  $w_-=0$ in a.e. $ O_-$, which is impossible with the definition of $ O_-$.\hfill$\Box$\smallskip

 We remark that when $L(s)=s^p$,  Lemma \ref{lm cp} could be seen in  \cite[Lemma 2.1]{CC}.
The following lemma plays an important role in the obtention of uniqueness for classical solution.

\begin{lemma}\label{cr hp}
Assume that $\Omega$ is a bounded $C^2$ domain and $L$ is a continuous function stated in Lemma \ref{lm cp} and satisfying
$L(x, 0)=0, \;\forall x \in \Omega$, then the homogeneous problem
\begin{equation}\label{eq0 2.1}
 \arraycolsep=1pt\left\{
\begin{array}{lll}
 \displaystyle \mathcal{L}_\mu u  +L(x,u)= 0\qquad
   {\rm in}\quad  {\Omega}\setminus \{0\},\\[1.5mm]
 \phantom{ L_\mu  +L(x,u)   }
 \displaystyle  u= 0\qquad  {\rm   on}\quad \partial{\Omega},\\[1.5mm]
 \phantom{   }
  \displaystyle \lim_{x\to0}u(x)\Phi_\mu^{-1}(x)=0
 \end{array}\right.
\end{equation}
has only zero solution.

\end{lemma}
{\bf Proof.}  Let $u$ be a solution of (\ref{eq0 2.1}), then for any $\epsilon>0$, there exists $r_\epsilon>0$ converging to zero as $\epsilon\to0$ such that
 $$u\le \epsilon \Phi_\mu\quad{\rm in}\quad \overline{B_{r_\epsilon}(0)}\setminus\{0\}.$$
We see that
$$u=0<\epsilon \Phi_\mu \quad{\rm on}\quad \partial\Omega,$$
then by Lemma \ref{lm cp}, we have that
$$u\le \epsilon \Phi_\mu\quad{\rm in}\quad \Omega\setminus\{0\}. $$
By the arbitrary of $\epsilon$, we have that $u\le 0$   in $\Omega\setminus\{0\}.$
The same way to obtain that $u\ge 0$   in $\Omega\setminus\{0\}.$
This ends the proof.\hfill$\Box$\medskip

 Now let
 $\mathcal{L}^*_\mu$ be defined by (\ref{L}).
 %=-\Delta -\frac{2 \tau_+(\mu) }{|x|^2}\,x\cdot\nabla $$ be
 As a consequence of the above result, an important test function we use often in the $d\mu$ distributional sense is the solution $\xi_0$ of
\begin{equation}\label{eq 2.2}
\arraycolsep=1pt\left\{
\begin{array}{lll}
 \displaystyle   \mathcal{L}_{\mu}^* u = 1\qquad
   {\rm in}\quad   \Omega,\\[2mm]
 \phantom{ \mathcal{L}_{\mu}^* }
 \displaystyle  u= 0\qquad  {\rm   on}\quad \partial{\Omega}
 \end{array}\right.
\end{equation}
satisfies
\begin{equation}\label{3.5}
0< \xi_0\le c\quad{\rm in} \quad  \Omega.
\end{equation}
In fact, without loss of the generality, let $\Omega\subset B_1(0)$ and denote
$v(r)=1-r^{2}$, then
$$ \mathcal{L}_{\mu}^* v=2N+2\tau_+(\mu)>1. $$
So $v$ is a super solution of (\ref{eq 2.2}).    For any $x_0\in\Omega$, since $\Omega$ is $C^2$ domain, there exists $r_0>0$ such that $B_{r_0}(x_0)\subset \Omega$.
Let $w_t(x)=t(r_0^2-|x-x_0|^2)$, then for some $t>0$ small,
$w_t$ is a sub solution of (\ref{eq 2.2}). Then  (\ref{3.5}) follows by  Lemma \ref{lm cp}.

Next, we build the distributional identity for the classical solution of the nonhomogeneous problem with "removable" singularity, i.e.  $\lim_{x\to0}u(x)\Phi_\mu^{-1}(x)=0$.

\begin{lemma}\label{lm 2.1}
Assume that  $f\in C^\gamma(\bar\Omega)$ for some $\gamma\in(0,1)$, then
\begin{equation}\label{2.02}
 \arraycolsep=1pt\left\{
\begin{array}{lll}
 \displaystyle  \mathcal{L}_\mu u= f\qquad
   {\rm in}\quad  {\Omega}\setminus \{0\},\\[1mm]
 \phantom{  L_\mu \, }
 \displaystyle  u= 0\qquad  {\rm   on}\quad \partial{\Omega},\\[1mm]
 \phantom{   }
  \displaystyle \lim_{x\to0}u(x)\Phi_\mu^{-1}(x)=0
 \end{array}\right.
\end{equation}
has a unique solution $u_f$
satisfying the distributional identity:
\begin{equation}\label{2.01}
 \int_{\Omega} u_f  \mathcal{L}_\mu^*(\xi )\, d\mu =\int_{\Omega} f \xi \, d\mu,\quad \forall\,\xi\in  C^{1.1}_0(\Omega).
\end{equation}
\end{lemma}
{\bf Proof.}  The uniqueness follows by Lemma \ref{cr hp} and  the existence could be derived by the Perron's method with the super solution $\bar u$  chosen as
$s(\Gamma_\mu(x)-t|x|^2)$ when $\mu_0\le \mu<2N$, $s(\Gamma_\mu(x)-t|x|^3)$ when $\mu=2N$
and $s|x|^2$ when $\mu>2N$, where $t>0$ is chosen such that the function is positive on $\Omega$ and then fix $s>0$ compared with $\norm{f}_{L^\infty(\Omega)}$, and the sub solution
$\underline{u}=-\bar u$.
We are going to prove the distributional identity.

\emph{The case $\mu>\mu_0$.} Indeed, for $\mu>\mu_0$, we can choose $\tau_0\in(\tau_-(\mu),\,\min\{2,\tau_+(\mu)\})$,
and denote
$$
V_0(x)=  |x|^{\tau_0},\quad\forall\, x\in \Omega\setminus\{0\}
$$
and
$$\mathcal{L}_\mu V_0(x)=c_{\tau_0} |x|^{\tau_0-2},$$
where
$$
c_{\tau_0} =\mu-\tau_0(\tau_0+N-2)>0.
$$
Since $f$ is bounded, there exists $t_0>0$ such that
$$|f(x)|\le t_0c_{\tau_0} |x|^{\tau_0-2},\qquad\forall x\,\in \Omega\setminus\{0\},$$
then $t_0 V_0$ and $-t_0V_0$ are super solution and sub solution of (\ref{2.02}), respectively.

For $n$ large,   $\Omega\setminus \overline{B_{\frac1n}(0)}$ is nonempty and the problem
\begin{equation}\label{eq2.2}
 \arraycolsep=1pt\left\{
\begin{array}{lll}
 \displaystyle  \mathcal{L}_\mu u= f \quad
   &{\rm in}\quad  \Omega\setminus \overline{B_{\frac1n}(0)},\\[2mm]
 \phantom{  - }
 u=0\quad  &{\rm   on}\quad \partial (\Omega\setminus \overline{B_{\frac1n}(0)})
 \end{array}\right.
\end{equation}
 has a  unique classical solution $w_n$.
 By Lemma \ref{lm cp}, we have that
$$|w_n(x)|\le t_0V_0(x),\quad\forall\, x\in\Omega\setminus\{0\}.$$
By standard arguments and regularity theory, taking $u_\mu=\lim_{n\to+\infty} w_n$, then $u_\mu$ is a classical solution of (\ref{2.02})
and
$$
| u_\mu(x)|\le t_0V_0(x),\quad\forall\, x\in\Omega\setminus\{0\}.
$$
From Lemma 4.9 in \cite{C} with $h\equiv1$ (e.g. \cite{T}), we have that
$$|\nabla u_\mu(x)|\le cV_0(x)|x|^{-1},\quad\forall\, x\in\Omega\setminus\{0\}.$$
Thus, for $\xi\in C^{1.1}_0(\Omega)$, multiplying $\Gamma_\mu\xi$ in (\ref{2.02}) and integrating over $\Omega\setminus \overline{B_r(0)}$,   we have that
\begin{eqnarray}
\int_{\Omega\setminus \overline{B_r(0)}} f\xi \, d\mu &=& \int_{\Omega\setminus \overline{B_r(0)}} \mathcal{L}_\mu(u_\mu) \xi   d\mu
    =   \int_{\Omega\setminus \overline{B_r(0)}}u_\mu\mathcal{ L}_\mu^*(\xi ) d\mu \nonumber \\&&
   +\int_{\partial B_r(0)}\left(\nabla u_\mu\cdot\frac{x}{|x|}  \Gamma_\mu  - \nabla  \Gamma_\mu  \cdot\frac{x}{|x|}u_\mu \right)\xi\,d\omega
 \label{02.1}  \\&& - \int_{\partial B_r(0)} u_\mu \Gamma_\mu  \left( \nabla\xi\cdot\frac{x}{|x|}\right) \,d\omega. \nonumber
\end{eqnarray}
For   $r=|x|>0$ small,
$$\left|\nabla u_\mu\cdot\frac{x}{|x|} \Gamma_\mu\,\xi  \right|\le c\norm{\xi }_{L^\infty(\Omega)} r^{\tau_0-1+\tau_+(\mu)} \quad{\rm
and}\quad  \left|\nabla \Gamma_\mu \cdot\frac{x}{|x|}u_\mu\,\xi\right|\le c\norm{\xi }_{L^\infty(\Omega)} r^{\tau_0+\tau_+(\mu)-1}.  $$
Since $\tau_0+\tau_+(\mu)-1>1-N$,
\begin{eqnarray*}
&&\lim_{r\to0^+}\int_{\partial B_r(0)} \left( \nabla u_\mu\cdot\frac{x}{|x|}  \Gamma_\mu  - \nabla  \Gamma_\mu  \cdot\frac{x}{|x|}u_\mu\right)\xi\,d\omega=0.
\end{eqnarray*}
For the last term in the right hand side of (\ref{02.1}), we have that
$$
\left|\int_{\partial B_r(0)} u_\mu \Gamma_\mu  \left( \nabla\xi\cdot\frac{x}{|x|}\right) \,d\omega \right|\le c_\mu \norm{\xi}_{C^1(\R^N)}r\to0\quad{\rm as}\quad r\to0^+.
$$
Therefore, we have that
\begin{equation}\label{2.03}
  \int_{\Omega}u_\mu\mathcal{L}^*_\mu(\xi)\,  d \mu = \int_{\Omega}f \xi\,  d \mu.
\end{equation}

\emph{The case $\mu=\mu_0$ and $\mu_0<0$.} %For $\mu=\mu_0$, we have to employee new method to prove the distributional identity.
 By the linearity of $\mathcal{L}_{\mu_0}$, we may assume that $f\in C^1(\bar\Omega)$ is nonnegative. By the Comparison Principle, we have that $u_\mu\ge 0$ in $\Omega\setminus\{0\}$.

We claim that $\mu\mapsto u_\mu$ is decreasing in $[\mu_0,0)$. Let $\mu_1\ge \mu_2>\mu_0$,  then
$$f=-\Delta u_{\mu_1} +\frac{\mu_1}{|x|^2}u_{\mu_1} \ge -\Delta u_{\mu_1} +\frac{\mu_2}{|x|^2}u_{\mu_1}, $$
so $u_{\mu_1}$ is a sub solution of (\ref{2.02}) with $\mu=\mu_2$. The monotonicity follows by the Comparison Principle in Lemma \ref{cr hp}.

We next construct a uniformly bound for $u_{\mu}$ for $\mu>\mu_0$.  Let
$$
V(x)=   |x| ^{ \tau_+(\mu_0)}-(s_0|x|)^2,\quad\forall x\in \Omega\setminus\{0\},
$$
where   $s_0>0$ and $V>0$ in $\Omega\setminus\{0\}$.
We see that
 $$\mathcal{L}_{\mu_0}V(x)=\frac{(N+2)^2}{4} s_0^2 >0,\quad\forall\, x\in\Omega\setminus\{0\}.$$
Then there exists $t_0>0$ such that
$$\mathcal{L}_{\mu_0}(t_0V)\ge \norm{f}_{L^\infty(\Omega)}. $$
By the Comparison Principle, we have that
$$u_\mu\le t_0V\quad {\rm in}\quad \Omega\setminus\{0\}.$$
For $\xi \in  C^{1.1}_0(\Omega)$, there exists $c>0$ independent of $\mu$ such that
$$|\mathcal{L}_\mu^*(\xi)|\le c\|\xi\|_{C^{1.1}_0(\Omega)}+ |\mu|\|\xi\|_{C^{1}_0(\Omega)} |x|^{-1}.%\frac{\morm{\xi}_{C^{1.1}_0(\Omega)}}{|x|}
$$

From the Dominate Monotonicity Convergence Theorem,  there exists $u_{\mu_0}\le t_0V$ such that
$$u_{\mu} \to u_{\mu_0} \quad{\rm as}\quad \mu\to \mu^+_0 \quad{\rm a.e.\ in}\ \ \Omega\quad {\rm and\ in}\quad L^1(\Omega,\, |x|^{-1}d\mu).$$
 Passing to the limit of (\ref{2.03}) as $\mu\to \mu^+_0$, we obtain that
 %\begin{equation}\label{2.04}
$$  \int_{\Omega}u_{\mu_0}  \mathcal{L}^*_{\mu_0}(\xi)\,  d\mu_0= \int_{\Omega}f \,\xi  d\mu_0,$$
% \end{equation}
%We finally prove the identity (\ref{2.04}) holds for any $\xi\in C^{1.1}_0(\Omega)$.  In fact, for $\xi\in C^{1.1}_0(\Omega)$ fixed,   let $\xi_n$ the solution of
%\begin{equation}\label{eq 2.1}
%\arraycolsep=1pt\left\{
%\begin{array}{lll}
% \displaystyle   \mathcal{L}_{\mu}^* u = f_n^+\qquad
%   {\rm in}\quad   \Omega,\\[2mm]
% \phantom{ \mathcal{L}_{\mu}^* }
% \displaystyle  u= 0\qquad  {\rm   on}\quad \partial{\Omega}.
% \end{array}\right.
%\end{equation}
%where  $$f_n^+=$$
which ends the proof.\hfill$\Box$

%\begin{remark}\label{re 2.1}
%In the special case that $\mu=\mu_0<0$, $f\equiv1$ and $\Omega=B_1(0)$, the solution of of (\ref{2.02}) is
%$\frac1{c_0}(\Gamma_{\mu_0}(x)-|x|^2)$, which is in the Hilbert space $H_{\mu_0}(\Omega)$, but it is not in normal Hilbert space $H^1_0(\Omega)$, the space of the closure of $C^{1.1}_0$ under the
%norm
%$$\norm{u}_{0} =\sqrt{\int_\Omega |\nabla u|^2 dx}.$$
%A direct consequence is that $H_{\mu_0}$ is not equivalent to $H^1_0(\Omega)$.

%\end{remark}

\subsection{Kato's inequality}

The following proposition is  the Kato's type estimate.

\begin{proposition}\label{pr 2.1}
Let  $\rho(x)={\rm dist}(x,\partial\Omega)$ and $f\in L^1(\Omega,\,\rho d\mu)$,  there exists a  unique $d\mu$-distributional
solution $u\in L^1(\Omega, |x|^{-1}d\mu)$ of the problem
\begin{equation}\label{homo}
 \arraycolsep=1pt\left\{
\begin{array}{lll}
  \mathcal{L}_\mu &u=f\quad & \rm{in}\quad\Omega,\\[2mm]
&u=0\quad & \rm{on}\quad \partial\Omega,
\end{array}\right.
\end{equation}
that is,
$$\int_\Omega u \mathcal{L}^*_\mu(\xi)\, d\mu =\int_\Omega f \xi\, d\mu,\quad\forall \, \xi\in C^{1.1}_0(\Omega). $$
Then for any $\xi\in  C^{1.1}_0(\Omega)$, $\xi\ge0$, we have that
 \begin{equation}\label{sign}
\int_\Omega |u|  \mathcal{L}_\mu^*(\xi)\, d\mu \le \int_\Omega
{\rm{sign}}(u)f  \xi\, d\mu
\end{equation}
and
 \begin{equation}\label{sign+}
\int_\Omega u_+  \mathcal{L}_\mu^*(\xi)\,  d\mu \le \int_\Omega
{\rm{sign}}_+(u)f \xi\, d\mu.
\end{equation}
%and
%\begin{equation}\label{sign0}
%\int_\Omega u_+(-\Delta)^\alpha \xi dx\le \int_\Omega  \xi
%\rm{sign}_+(u) \nu dx,
%\end{equation}
%where $u_+=\max\{0,u\}$, $\rm{sign}_+(u)=1$ if $u>0$ and
%$\rm{sign}_+(u)=0$ if $u\le0$.
\end{proposition}
{\bf Proof.} {\it Uniqueness}.
Let $w$ be a $d\mu$-distributional solution of
\begin{equation}\label{L0}\left\{\begin{array}{lll}
\mathcal{L}_\mu w=0\qquad\mbox{ in }\quad\Omega,\\
\phantom{ \mathcal{L}_\mu  }w=0\qquad\mbox{ on }\quad\partial\Omega.
\end{array}\right.\end{equation}
For any Borel subset $O$  of $\Omega$, denote by $\eta_{\omega,n}$ the
solution of
\begin{equation}\label{L00}
\arraycolsep=1pt\left\{
\begin{array}{lll}
 \displaystyle  \mathcal{L}_\mu^*  u= \zeta_n\qquad
   &{\rm in}\quad   \Omega,\\[1mm]
 \phantom{  \mathcal{L}_\mu^*  }
 \displaystyle  u= 0\quad & {\rm   on}\quad \partial{\Omega},
 \end{array}\right.
 \end{equation}
where $\zeta_n:\bar\Omega\mapsto[0,1]$ is a $C^1(\bar\Omega)$
function such that $\zeta_n\to\chi_O\; {\rm{in}}\ L^\infty( \Omega)\; {\rm{as}}\ n\to\infty.$
By Lemma \ref{lm 2.1}, we have that $\displaystyle\int_{\Omega}w\Gamma_\mu\, \zeta_n \ dx=0.
$
Passing to the limit as $n\to\infty$, we have that
$$\displaystyle\int_{O}w\Gamma_\mu \ dx=0,$$
which implies that $w=0$.
\smallskip

\noindent{\it Existence and estimate (\ref{sign})}. For  $\sigma>0$,
we define an even convex function  $\phi_\sigma$ as
\begin{equation}\label{te2}
\phi_\sigma(t)=\left\{ \arraycolsep=1pt
\begin{array}{lll}
|t|-\frac\sigma2\quad & \rm{if}\quad |t|\ge \sigma ,\\[2mm]
\frac{t^2}{2\sigma}\quad & \rm{if} \quad |t|< \sigma/2.
\end{array}
\right.
\end{equation}
 Then for any $t\in\R$,  $|\phi_\sigma'(t)|\le1$, $\phi_\sigma(t)\to|t|$ and
$\phi_\sigma'(t)\to\rm{sign}(t)$ when $\sigma\to0^+$.
%%%%%%%%%%%%%%%%%%%%%%%%%%%%%%%%%%%%%%%%%%%%%%%%%%%%%%%%%%%%%%%%%%%%%%%%%%%%%%%%%%%%%%%%%%%%%%%%%%%%%%%%%%%%%%%%%%%%%%%%%%%%%%%%%%%%%%%%%%%%%%
Let $\{f_n\}$ be a sequence of functions in $C^1(\bar\Omega)$ such
that
$$\lim_{n\to\infty} \int_\Omega|f_n-f|\rho  d\mu=0.$$
Let $u_n$ be the corresponding solution  to (\ref{homo}) with
right-hand side $f_n$, then   for any $\sigma>0$ and $\xi\in    C^{1.1}_0(\Omega),\ \xi\ge0$, we have that
\begin{equation}\label{2.7}
 \arraycolsep=1pt
\begin{array}{lll}
\displaystyle\int_\Omega \phi_\sigma(u_n) (-\Delta) (\Gamma_\mu\xi)\,
dx&=\displaystyle\int_\Omega  \Gamma_\mu \xi(-\Delta) \phi_\sigma(u_n)\,dx\\[4mm]
 &=\displaystyle  \int_\Omega
   \xi\phi_\sigma'(u_n)  (-\Delta)  u_n d\mu- \int_\Omega
   \xi\phi_\sigma''(u_n)  |\nabla u_n|^2 d\mu
\\[4mm]\displaystyle
&\le \displaystyle   \int_\Omega  \xi \phi_\sigma'(u_n) f_nd\mu-\mu\int_\Omega   \phi_\sigma'(u_n)\frac{ u_n}{|x|^2} \xi\, d\mu.
\end{array}
\end{equation}
Letting $\sigma\to 0$, we obtain that
$$\int_\Omega |u_n|(-\Delta) (\Gamma_\mu \xi) dx\leq  \int_\Omega   \xi \ {\rm {sign}} (u_n) f_n \ d\mu-\mu \int_\Omega   \frac{ |u_n|}{|x|^2}  \xi\, d\mu,$$
that is,
\begin{equation}\label{L1}
 \int_\Omega |u_n| \mathcal{L}_\mu^* (\xi) d\mu\leq  \int_\Omega  \xi |f_n|\, d\mu.
\end{equation}
Let $\eta_1$ be the solution of (\ref{eq 2.2}), taking $\xi=\xi_0$, the solution of (\ref{eq 2.2}),  then
\begin{equation}\label{L2}
\displaystyle\int_\Omega |u_n|\,  d\mu \leq  c_8\displaystyle \int_\Omega
|f_n|  \rho\, d\mu.
\end{equation}
Similarly,
\begin{equation}\label{L3}
\displaystyle\int_\Omega |u_n-u_m|   d\mu\leq  c_8\displaystyle
\int_\Omega |f_n-f_m|  \rho\, d\mu.
\end{equation}
Therefore, $\{u_n\}$ is a Cauchy sequence in $L^1(\Omega,\,   \rho d\mu )$ and its limit $u$
is  a $d\mu$-distributional solution of (\ref{homo}). Passing to the limit as $n\to\infty$ in
$(\ref{L1})$, we obtain (\ref{sign}).   Inequality (\ref{sign+}) is
proved by replacing $\phi_\sigma$ by $\tilde\phi_\sigma$ which is
zero on $(-\infty,0)$ and $\phi_\sigma$ on $[0,\infty)$.
\hfill$\Box$\medskip

\setcounter{equation}{0}
\section{Distributional identities}

\subsection{ Fundamental solution.}
The distributional identity of the fundamental solutions is derived by divergence theorem. \smallskip

 \noindent {\bf Proof of Theorem \ref{teo 1}.} For $\xi\in C^{1.1}_c(\R^N)$, multiply $\Gamma_\mu\xi$ in (\ref{eq 1.1}) and integrate over $\R^N\setminus \overline{B_r(0)}$,  we have that
\begin{eqnarray*}
0 &=& \int_{\R^N\setminus \overline{B_r(0)}} \mathcal{L}_\mu (\Phi_\mu) \Gamma_\mu \xi\, dx
    =   \int_{\R^N\setminus \overline{B_r(0)}}\Phi_\mu \mathcal{L}_\mu^*(\xi )\, d\mu  \\&&
   +\int_{\partial B_r(0)}\left( \nabla \Phi_\mu\cdot\frac{x}{|x|}  \Gamma_\mu  -\nabla  \Gamma_\mu  \cdot\frac{x}{|x|}\Phi_\mu\right)\xi\,d\omega
- \int_{\partial B_r(0)} \Phi_\mu \Gamma_\mu  \left( \nabla\xi\cdot\frac{x}{|x|}\right) \,d\omega.
\end{eqnarray*}
For $\mu>\mu_0$, we see that
\begin{eqnarray*}
 \nabla \Phi_\mu(x)\cdot\frac{x}{|x|} \Gamma_\mu(x)- \nabla \Gamma_\mu(x) \cdot\frac{x}{|x|}\Phi_\mu(x) %&=& (\tau_-(\mu)-\tau_+(\mu)) |x|^{\tau_-(\mu)+\tau_+(\mu)-1} \\
     &=&  - 2\sqrt{\mu-\mu_0}\, |x|^{1-N}
\end{eqnarray*}
and for $r=|x|>0$ small,
 $$|\xi(x)-\xi(0)|\le \norm{\nabla\xi }_{L^\infty(\R^N)} r$$
and
\begin{eqnarray*}
 \int_{\partial B_r(0)}\left( \nabla \Phi_\mu\cdot\frac{x}{|x|}  \Gamma_\mu  - \nabla  \Gamma_\mu  \cdot\frac{x}{|x|}\Phi_\mu\right)\xi(0)\,d\omega = - c_{\mu} \xi(0),
\end{eqnarray*}
where $c_\mu = 2\sqrt{\mu-\mu_0} |{S}^{N-1}|$ for $\mu > \mu_0$.%, and
%\begin{eqnarray*}
% \left| \int_{\partial B_r(0)} \left( \nabla \Phi_\mu\cdot\frac{x}{|x|}  \Gamma_\mu  - \nabla  \Gamma_\mu  \cdot\frac{x}{|x|}\Phi_\mu\right)\norm{ \xi }_{C^1(\R^N)} r\,d\omega \right|
%   \le c_{\mu} \norm{\nabla\xi }_{L^\infty(\R^N)} r,
%\end{eqnarray*}
%then,
%$$
% \lim_{r\to0}\left(\int_{\partial B_r(0)} \nabla \Phi_\mu\cdot\frac{x}{|x|}  \Gamma_\mu\xi \,d\omega
%   - \int_{\partial B_r(0)} \nabla  \Gamma_\mu \cdot\frac{x}{|x|}\Phi_\mu\xi \,d\omega\right)=- c_\mu\xi(0).
%$$
%Moreover,
%$$ \left|\int_{\partial B_r(0)} \Phi_\mu \Gamma_\mu  \left( \nabla\xi\cdot\frac{x}{|x|}\right) \,d\omega \right|\le c_\mu \norm{\nabla\xi }_{L^\infty(\R^N)} r. $$

Elementary estimates show  that
$$\lim_{r\to0^+}\int_{\R^N\setminus \overline{B_r(0)}}\Phi_\mu \mathcal{L}_\mu^*( \xi ) d\mu=c_\mu  \xi(0).$$

For $\mu=\mu_0$,  we have that $\tau(\mu_0)=\bar\tau(\mu_0)=\frac{2-N}{2}$,
 \begin{eqnarray*}
\nabla \Phi_{\mu_0}\cdot\frac{x}{|x|}\,  \Gamma_{\mu_0}  - \nabla \Gamma_{\mu_0} \cdot\frac{x}{|x|}\,\Phi_{\mu_0}=-|x|^{1-N},
 \end{eqnarray*}
 then
  $$\lim_{r\to0}\int_{\partial B_r(0)}\left( \nabla \Phi_{\mu_0}\cdot\frac{x}{|x|} \Gamma_{\mu_0}\xi
   -   \nabla  \Gamma_{\mu_0} \cdot\frac{x}{|x|}\Phi_{\mu_0}\xi \right) d\omega=  |\mathcal{S}^{N-1}|\xi(0)$$
 and
$$\lim_{r\to0} \int_{\partial B_r(0)} \Phi_{\mu_0} \Gamma_{\mu_0}  \left( \nabla\xi\cdot\frac{x}{|x|}\right) \,d\omega =0.$$%$$\left|\int_{\partial B_r(0)} \Phi_{\mu_0} \Gamma_{\mu_0}  \left( \nabla\xi\cdot\frac{x}{|x|}\right) \,d\omega\right|\le c_{\mu_0} \norm{\nabla\xi }_{L^\infty(\R^N)} (-\ln r)r,$$
%with $c_{\mu_0}=|\mathcal{S}^{N-1}|$.
Thus, we have that
$$\lim_{r\to0^+}\int_{\R^N\setminus \overline{B_r(0)}}\Phi_{\mu_0} \mathcal{L}_{\mu_0}^*(\xi)  d{\mu_0}=c_{\mu_0}  \xi(0).$$
This ends the proof.\hfill$\Box$\medskip

 \begin{remark}
For $\mu\not=0$, it is not proper to divide into each term like
$$\lim_{r\to0^+} \int_{ \R^N \setminus \overline{B_r(0)}}\Phi_\mu \mathcal{L}_\mu(\Gamma_\mu\xi) dx=\lim_{r\to0^+}\int_{ \R^N \setminus \overline{B_r(0)}}\Phi_\mu(-\Delta)(\Gamma_\mu\xi) dx
+ \lim_{r\to0^+}\int_{\R^N \setminus \overline{B_r(0)}}\mu \frac{\Phi_\mu  \Gamma_\mu\xi}{ |x|^{2}} dx.$$
In fact, we see that if $\xi(0)>0$,
$$\lim_{r\to0^+}\int_{ \R^N \setminus \overline{B_r(0)}}\frac{\Phi_\mu  \Gamma_\mu\xi}{ |x|^{2}} dx=+\infty.$$
\end{remark}

\subsection{$d\mu$-distributional solution in bounded domain}

Different from the case in $\R^N$, there is only one branch of  the fundamental solution for $\mathcal{L}_\mu u=\delta_0$, subject to the Dirichlet boundary condition.\medskip

\noindent {\bf Proof of Theorem \ref{teo 2}.} The existence follows from Perron's methods with the super and sub-solutions
$\Phi_\mu$ and  $\Phi_\mu-t\Gamma_\mu$ respectively. Here $t$ is taken such that $\Phi_\mu-t\Gamma_\mu\le 0$ on $\partial\Omega$. The uniqueness follows from the
Comparison Principle.

  Let $\eta_0:[0,+\infty)\to [0,\,1]$ be a decreasing $C^\infty$ function such that
\begin{equation}\label{eta}
 \eta_0=1\quad{\rm in}\quad [0,1]\qquad{\rm and}\qquad \eta_0=0\quad {\rm in}\quad[2,+\infty).
\end{equation}
Take $n_0\ge 1$ such that
$$\frac1{n_0}\sup\{r>0:\, B_r(0)\subset \Omega\}\le \frac12.$$
Denote $\eta_{n_0}(r)=\eta_0(n_0r)$ for $r\ge0$ and  $w_1=\Phi_\mu \eta_{n_0}$ and  $w_2=G_\mu-\Phi_\mu \eta_{n_0}$, then
\begin{equation}\label{2.3}
\arraycolsep=1pt\left\{
\begin{array}{lll}
 \displaystyle  \mathcal{L}_\mu w_i= -\nabla \eta_{n_0}\cdot \nabla  \Phi_\mu-\Phi_\mu \Delta\eta_{n_0}\qquad
   {\rm in}\quad  {\Omega}\setminus \{0\},\\[2mm]
 \phantom{  L_\mu \, }
 \displaystyle  w_i= 0\qquad  {\rm   on}\quad \partial{\Omega},\\[2mm]
 \phantom{   }
  \displaystyle \lim_{x\to0}w_i(x)\Phi_\mu^{-1}(x)= 2-i,
 \end{array}\right.
\end{equation}
where $i=1,2$.
We see that $-\nabla \eta_{n_0}\cdot \nabla  \Phi_\mu-\Phi_\mu \Delta\eta_{n_0}$ has compact set in $\overline{B_{\frac2{n_0}}(0)}\setminus B_{\frac1{n_0}}(0)$
and then  $-\nabla \eta_{n_0}\cdot \nabla  \Phi_\mu-\Phi_\mu \Delta\eta_{n_0}$ is smooth and bounded.

For $i=1$, following the proof of Theorem \ref{teo 1}, we have that
\begin{equation}\label{2.4}
 \int_{\R^N}w_1 \mathcal{L}_\mu^*( \xi)\, d\mu =\int_{\Omega} \left[-\nabla \eta_{n_0}\cdot \nabla  \Phi_\mu-\Phi_\mu \Delta\eta_{n_0}\right] \xi\, d\mu+c_\mu \xi(0),\quad \forall\,\xi\in C^{1.1}_0(\Omega).
\end{equation}
For $i=2$,  it follows by Lemma \ref{lm 2.1}  that
\begin{equation}\label{2.5}
 \int_{\Omega} w_2 \mathcal{L}_\mu^*( \xi ) d\mu=\int_{\Omega} \left[-\nabla \eta_{n_0}\cdot \nabla  \Phi_\mu-\Phi_\mu \Delta\eta_{n_0}\right] \xi\,d\mu,\quad \forall\,\xi\in C^{1.1}_0(\Omega).
\end{equation}
Together with (\ref{2.4}),
$$ \int_{\Omega}G_\mu \mathcal{L}_\mu^*( \xi)\, d\mu =c_\mu \xi(0),\quad \forall\,\xi\in  C^{1.1}_0(\Omega).$$
This completes the proof.\hfill$\Box$

\subsection{Approximation of the fundamental solution}

In this subsection, we see the approximation of the fundamental solution $G_\mu$.

\begin{proposition}\label{lm 2.3}
Let $\{\delta_n\}_n$ be a sequence  of nonnegative  $L^\infty$-functions such that supp$\,\delta_n\subset B_{r_n}(0)$, where $r_n\to0$ as $n\to+\infty$ and
$$ \delta_n \to \delta_0\quad{\rm as} \quad n\to+\infty\ \ {\rm in\ the\ distributional\ sense.}$$
 For any $n$, let $w_n$  be the unique solution of the problem in the $d\mu$-distributional sense
 \begin{equation}\label{eq 2.3}
 \arraycolsep=1pt\left\{
\begin{array}{lll}
 \displaystyle  \mathcal{L}_\mu u = \delta_n\qquad
   {\rm in}\quad  {\Omega}\setminus \{0\},\\[1mm]
 \phantom{  L_\mu -- }
 \displaystyle  u= 0\qquad  {\rm   on}\quad \partial{\Omega},\\[1mm]
 \phantom{   }
  \displaystyle \lim_{x\to0}u(x)\Phi_\mu^{-1}(x)=0.
 \end{array}\right.
\end{equation}
Then
$$\lim_{n\to+\infty} w_n(x)=\frac1{c_\mu}G_\mu(x),\quad \forall \, x\in\Omega\setminus\{0\}$$
and for any compact set $K\subset \Omega\setminus\{0\}$,
 \begin{equation}\label{3.7}
  w_n\to \frac1{c_\mu} G_\mu \quad{\rm as}\quad n\to+\infty \ \ {\rm in}\ \ C^2(K).
 \end{equation}
\end{proposition}
{\bf Proof.} From Lemma \ref{lm 2.1}, equation (\ref{eq 2.3}) has unique solution $w_n\ge0$ satisfying
\begin{equation}\label{3.0}
  \int_{\Omega} w_n  \mathcal{L}_\mu^*(\xi )\, d\mu =\int_{\Omega} \delta_n\xi\,  dx,\quad \forall\,\xi\in C^{1.1}_0(\Omega).
\end{equation}
In particular, by taking $\xi=\xi_0$, the solution of (\ref{eq 2.2}), we obtain that
$$\norm{w_n}_{L^1(\Omega,\,  d\mu)}\le \norm{\xi_0}_{L^\infty(\Omega)}\norm{\delta_n}_{L^1(\Omega)}=\norm{\xi_0}_{L^\infty(\Omega)}. $$
For any  $r>0$, take $\xi$ with the support in $\Omega\setminus B_r(0)$,
then  $\xi\in C^{1.1}_c(\overline{\Omega\setminus B_r(0)})$,
$$\int_{\Omega\setminus B_r(0)} w_n  \mathcal{L}_{\mu}^* (\xi) \, d\mu =0. $$

We claim that  $w_n$ is uniform bounded in $L^1(\Omega, \,|x|^{-1} d\mu)$. Indeed, for  $\sigma>0$,  let
\begin{equation}\label{t1}
\xi_\sigma(x)= (r_0-\phi_\sigma(|x|))  \eta_{n_0}(|x|),\quad\forall x\in\Omega,
\end{equation}
where $\phi_\sigma$ is defined in (\ref{te2}) and $\eta_{0}$ is given in (\ref{eta}) and $n_0$ is such that $\eta_{n_0}$ has compact set in $\Omega$ and $r_0>0$ is such that
$\Omega\subset B_{r_0}(0)$.
Then $\xi_\sigma\in C^{1.1}_0(\Omega)$ and for $x\in B_{\frac1{n_0}}(0)$, we have that
 \begin{equation}\label{eq 3.3000}
 \mathcal{L}_\mu^*(\xi_\sigma )(x)=\arraycolsep=1pt\left\{
\begin{array}{lll}
 \frac{N-1+\tau_+(\mu)}{|x|}\quad&{\rm for} \quad  \sigma\le |x|<1/n_0,\\[2mm]
 \frac1\sigma+N-1+\tau_+(\mu)\quad&{\rm for} \quad  |x|<\sigma
\end{array}
\right.
\end{equation}
and
there exists $c>0$ such that
 \begin{equation}\label{eq 3.3100}
|\mathcal{L}_\mu^*(\xi_\sigma )(x)|\le c,\qquad\forall\, x\in \Omega \setminus B_{\frac1{n_0}}(0),
\end{equation}
where $N-1+\tau_+(\mu)>0$ for $\mu\ge\mu_0$.
Therefore, we have that
\begin{eqnarray*}
(N-1+\tau_+(\mu))\int_{B_{\frac1{n_0}}(0)\setminus B_\sigma(0)} w_n\,  |x|^{-1} d\mu &\le&
 \int_{B_{\frac1{n_0}}(0)} w_n\,  \mathcal{L}_\mu^*(\xi_\sigma ) d\mu
 \\ &\le & \int_\Omega \delta_n \xi_\sigma  d\mu +c\int_{\Omega\setminus B_{\frac1{n_0}}(0)}w_n   d\mu
 \\ &\le & c,
\end{eqnarray*}
where $c>0$ independent of $n$.
Passing to the limit as $\sigma\to0$, we have that
$$\int_{B_{\frac1{n_0}}(0)} w_n\,  |x|^{-1} d\mu\le c',$$
which, together with the fact that
$$\int_{\Omega\setminus B_{\frac1{n_0}}(0)} w_n \,  |x|^{-1} d\mu\le  c,$$
implies that
$$\int_{\Omega} w_n\,  |x|^{-1} d\mu\le c',$$
that is,  $w_n$ is uniform bounded in $L^1(\Omega, \,|x|^{-1} d\mu)$.

From Corollary 2.8 in \cite{V} with $L^*=\mathcal{L}_{\mu}^*$, which is strictly elliptic in $\Omega\setminus B_r(0)$, we have that for $q<\frac{N}{N-1}$,
\begin{eqnarray*}
\norm{w_n\Gamma_\mu}_{W^{1,q}(\Omega_{2r})}  &\le & c\norm{\delta_n}_{L^1 (\Omega\setminus B_r(0))}+ c\norm{w_n }_{L^1(\Omega\setminus B_r(0),\,d\mu)} \\
  &\le & c(1+\norm{\xi_0}_{L^\infty(\Omega)}),
\end{eqnarray*}
where $\Omega_{2r}=\{x\in\Omega\setminus B_{2r}(0):\, \rho(x)>2r\}.$
By the compact embedding
$$W^{1,q}(\Omega_{2r})\hookrightarrow L^1(\Omega_{2r}),$$
up to some subsequence, there exists $w_\infty\in W^{1,q}_{loc}(\Omega)\cap L^1(\Omega,\, d\mu)$ such that
$$w_n \to w_\infty  \quad{\rm as}\quad n\to+\infty\quad {\rm a.e.\ \ in}\ \ \Omega\ \ {\rm and\ in}\quad L^1(\Omega,\,d\mu)$$
and it follows by (\ref{3.0}) that for $\xi\in C^{1.1}_0(\Omega)$,
$$\int_{\Omega} w_\infty   \mathcal{L}_\mu^*(\xi)\, d\mu =\xi(0). $$
Furthermore,
$$\int_{\Omega}( w_\infty-\frac1{c_\mu}G_\mu)   \mathcal{L}_\mu^*(\xi)\, d\mu =0. $$
From the Kato's inequality Lemma \ref{pr 2.1} with $f=0$, we deduce that
$$w_\infty=\frac1{c_\mu}G_\mu\quad{a.e.}\;\; \Omega. $$

{\it Proof of (\ref{3.7}).} For any $x_0\in \Omega\setminus\{0\}$, let $r_0=\frac14\{|x_0|,\, \rho(x_0)\}$ and
$$v_n=w_n\eta,$$
where $\eta(x)=\eta_0(\frac{|x-x_0|}{r_0})$. There exists $n_0>0$ such that for $n\ge n_0$,
$${\rm supp}v_n\cap B_{r_n}(0)=\emptyset.$$
Then
\begin{eqnarray*}
-\Delta v_n (x) &=& -\Delta w_n (x)\eta(x)-\nabla w_n\cdot\nabla\eta-w_n \Delta\eta \\
   &=&  -\nabla w_n\cdot\nabla\eta-w_n \Delta\eta,
\end{eqnarray*}
where $\nabla\eta$ and $\Delta\eta$ are smooth.

We observe that $w_n\in W^{1,q}(B_{2r_0}(x_0))$ and
$-\nabla w_n\cdot\nabla\eta-w_n \Delta\eta\in L^q(B_{2r_0}(x_0)),$
then we have that
$$\norm{v_n}_{ W^{2,q}(B_{r_0}(x_0))}\le c\norm{w_n}_{L^1(\Omega,\,  d\mu)},$$
where $c>0$ is independent of $n$.
Thus,  $-\nabla w_n\cdot\nabla\eta-w_n \Delta\eta\in W^{1,q}(B_{r_0}(x_0)),$
repeat above process $N_0$ steps, for $N_0$ large enough, we deduce that
$$\norm{w_n}_{C^{2,\gamma}(B_{\frac{r_0}{2^{N_0}}}(x_0))}\le c\norm{w_n}_{L^1(\Omega,\,  d\mu)},$$
where $\gamma\in(0,1)$ and $c>0$ is independent of $n$.
As a conclusion, (\ref{3.7}) follows by Arzel\`{a}-Ascola theorem and Heine-Borel theorem.
This ends the proof. \hfill$\Box$

\begin{remark}
Let $\eta_0:[0,+\infty)\to[0,1]$ be a $C^\infty$ decreasing  function such that $\eta_0=1$ in $[0,1]$ and $\eta_0=0$ in $[2,+\infty)$,  $\delta_n(x)=n^N\eta_0(n|x|)$ for $x\in\R^N$.  Then
$$\delta_n\to \delta_0\quad {\rm in} \;\;\mathcal D'(\R^N)\;\; {\rm as}\quad n\to+\infty.$$
%\quad{\rm in\ the\ distribution\ sense.

\end{remark}
\setcounter{equation}{0}
\section{Nonhomogeneous problem}

\subsection{Isolated singularities}
We concerns with the isolated singularities of the solution $v_f$ of (\ref{eq 1.1f}) verifying the identity (\ref{1.2f}) with $k=0$.
Let  $\mathcal{G}_\mu$ be the Green's kernel of $\mathcal{L}_\mu$ in $\Omega\times\Omega$ with Dirichlet boundary condition.
From Theorem \ref{teo 2},   it holds that for $x\in\Omega\setminus\{0\}$ and $y=0$,
$$\mathcal{G}_\mu(x,0)=G_\mu(x),$$
which expresses $\delta_0$ in the $d\mu$-distributional sense. For $x, y\in \Omega\setminus\{0\}$, we have following estimate of $\mathcal{G}_\mu$.

\begin{lemma}\label{lm 4.2.1}
If $\mu_0<0$ and $\mu_0\le \mu<0$, we have that for $x, y\in \Omega\setminus\{0\}$,
\begin{equation}\label{e g1}
 0<\mathcal{G}_\mu(x,y)\le %\arraycolsep=1pt\left\{
%\begin{array}{lll}
c\left(|x-y|^{2-N}+\frac{|x|^{\tau_+(\mu)} }{|x-y|^{N-2+\tau_+(\mu)}}+\frac{ |y|^{\tau_+(\mu)}}{|x-y|^{N-2+\tau_+(\mu)}}+ \frac{|x|^{\tau_+(\mu)}|y|^{\tau_+(\mu)}}{|x-y|^{N-2+2\tau_+(\mu)}}\right);%\ \, \forall x,y\in \Omega\setminus\{0\};\\[2mm]
%\Phi_\mu(x),\qquad x\in \Omega\setminus\{0\},\quad y=0;
%\end{array}
%\right.
\end{equation}
if $\mu\ge0$ and $N\ge3$, we have that for $x, y\in \Omega\setminus\{0\}$,
 \begin{equation}\label{e g2}
 0<\mathcal{G}_\mu(x,y)\le %\arraycolsep=1pt\left\{
%\begin{array}{lll}
c\min\left\{|x-y|^{2-N},\,\frac{|x|^{\tau_+(\mu)} }{|x-y|^{N-2+\tau_+(\mu)}}, \frac{ |y|^{\tau_+(\mu)}}{|x-y|^{N-2+\tau_+(\mu)}},  \frac{|x|^{\tau_+(\mu)}|y|^{\tau_+(\mu)}}{|x-y|^{N-2+2\tau_+(\mu)}}\right\}.%\\[2mm]
%\Phi_\mu(x),\qquad x\in \Omega\setminus\{0\},\quad y=0.
%\end{array} \right.
 \end{equation}
Here and in the following, $|x-y|^{2-N}$ in (\ref{e g2}) should be replaced by $-\ln|x-y|$ when $N=2$.

\end{lemma}

{\bf Proof.} For $\mu_0<0$ and $\mu_0\le \mu<0$, it follows from Theorem 3.11 in \cite{FMT} that the corresponding heat kernel verifies
\begin{eqnarray*}
k_\mu(t,x,y) &\le & c\min\left\{(1+\frac{|x|}{\sqrt{t}})^{\tau_+(\mu)}(1+\frac{|y|}{\sqrt{t}})^{\tau_+(\mu)}, (\frac{d(x)d(x)}{|x||y|})^{\tau_+(\mu)} \right\} t^{-\frac N2} e^{-c\frac{|x-y|^2}{t}} \\
    &\le &c\left(1+ (\frac{|x|}{\sqrt{t}})^{\tau_+(\mu)}+(\frac{|y|}{\sqrt{t}})^{\tau_+(\mu)}+(\frac{|x||y|}{t})^{\tau_+(\mu)}\right)t^{-\frac N2} e^{-c\frac{|x-y|^2}{t}},
\end{eqnarray*}
which, together with $\mathcal{G}(x,y)=\displaystyle\int_0^\infty k_\mu(t,x,y)dt$, implies (\ref{e g1}).

For $\mu \ge0$, it follows from Proposition 1.1 in \cite{IYM} that
\begin{eqnarray*}
k_\mu(t,x,y) &\le &c \min\left\{1,\,(\frac{|x|}{\sqrt{t}})^{\tau_+(\mu)}\}\, \min\{1,\,(\frac{|y|}{\sqrt{t}})^{\tau_+(\mu)}\right\}\,t^{-\frac N2} e^{-c\frac{|x-y|^2}{t}} \\
    &\le &c\min\left\{1,\, (\frac{|x|}{\sqrt{t}})^{\tau_+(\mu)},\,(\frac{|y|}{\sqrt{t}})^{\tau_+(\mu)},\,(\frac{|x||y|}{t})^{\tau_+(\mu)}\right\}\,t^{-\frac N2} e^{-c\frac{|x-y|^2}{t}},
\end{eqnarray*}
then  (\ref{e g2}) holds.\hfill$\Box$\medskip

\begin{remark}\label{re 4.1}  % \textcolor[rgb]{1.00,0.00,0.00}{check the case} $\mu=\mu_0$.
For $\mu_0<0$ and $\mu_0\le \mu<0$, it follows from \cite{FMT} that
the heat kernel has the lower bound as
$$k_\mu(t,x,y) \ge  c\min\left\{(1+\frac{|x|}{\sqrt{t}})^{\tau_+(\mu)}(1+\frac{|y|}{\sqrt{t}})^{\tau_+(\mu)}, (\frac{d(x)d(x)}{|x||y|})^{\tau_+(\mu)} \right\} t^{-\frac N2} e^{-c\frac{|x-y|^2}{t}},$$
then we have that the Green kernel has lower estimates: for any compact set $K$ in $\Omega$, there exists $c>0$ such that for $x,y\in K\setminus\{0\}$,
$$\mathcal{G}_\mu(x,y)\ge
c\left(|x-y|^{2-N}+\frac{|x|^{\tau_+(\mu)} }{|x-y|^{N-2+\tau_+(\mu)}}+\frac{ |y|^{\tau_+(\mu)}}{|x-y|^{N-2+\tau_+(\mu)}}+ \frac{|x|^{\tau_+(\mu)}|y|^{\tau_+(\mu)}}{|x-y|^{N-2+2\tau_+(\mu)}}\right).$$

\end{remark}

\begin{proposition}\label{pr 4.1.1}
Let    $f$ be a function in $C^1(\overline{\Omega}\setminus \{0\})$ satisfying (\ref{f1}) and (\ref{4.1-3}).

Denote by $\mathbb{G}_\mu$ the Green operator defined by
$$\mathbb{G}_\mu[f](x)=\int_\Omega \mathcal{G}(x,y)f(y) dy\quad{\rm for}\ \ f\in L^1(\Omega,\, d\mu).$$
Then  $\mathbb{G}_\mu[f]$ is the classical solution of (\ref{eq 1.1f}) verifying
\begin{equation}\label{1.01}
 \lim_{x\to0}\mathbb{G}_\mu[f](x)\Phi_\mu^{-1}(x)=0,
\end{equation}

\end{proposition}
\noindent{\bf Proof.} From the Green's kernel, it is known that $\mathbb{G}_\mu[f]$ is the solution of (\ref{eq 1.1f}). We next prove (\ref{1.01})
and we may assume that $f$ is nonnegative, if not, we just replace $f$ by $|f|$  in the following proof.

 {\it Case 1: $N\ge 3$ and $\mu\ge0$.} When $N=2$, there are just some small differences, we omit the proof.

   For $x\in\Omega\setminus\{0\}$,
 \begin{eqnarray}
\mathbb{G}_\mu [f](x)\Phi_\mu^{-1}(x)  &\le&  c \left( \Phi_\mu^{-1}(x)  \int_{B_{r_0}(0)\setminus B_{\frac{|x|}2}(x)}
\frac{|y|^{\tau_+(\mu)} f(y)}{|x-y|^{N-2+\tau_+(\mu)}}  dy+ \Phi_\mu^{-1}(x) \int_{B_{\frac{|x|}2}(x)}
\frac{ f(y)}{|x-y|^{N-2}}  dy\right)\nonumber
\\&=:& c \left(  \mathbb{E}_1(x)+\mathbb{E}_2(x)\right),\label{4.5}
 \end{eqnarray}
 where $r_0=\sup_{x,y\in\Omega}|x-y|$ and $f$ is extended to be that $f=0$ in $B_{r_0}(0)\setminus \Omega$.

 We observe that
 \begin{eqnarray*}
  \int_{B_{\frac{|x|}2}(x)} \frac{ f(y)}{|x-y|^{N-2}}  dy &\le& \sup_{z\in B_{\frac{3|x|}2}(0)\setminus B_{\frac{|x|}2}(0) } f(z)\, \int_{B_{\frac{|x|}2}(x)} |x-y|^{2-N}  dy
  \\&\le& c|x|^{2}\sup_{B_{\frac{3|x|}2}(0)\setminus B_{\frac{|x|}2}(0)} f(z)
\end{eqnarray*}
and
\begin{equation}\label{4.1-31}
 \mathbb{E}_2(x) \le  c  \sup_{z\in B_{\frac{3|x|}2}(0)\setminus B_{\frac{|x|}2}(0)}|z|^{2-\tau_-(\mu)} f(z) \to0\quad {\rm as}\quad |x|\to0
\end{equation}
by the assumption (\ref{4.1-3}).

For $y\in B_{r_0}(0)\setminus B_{\frac{|x|}2}(x)$, we have that
$$|x-y|\ge\frac18 (|x|+|y|)  $$ and
\begin{eqnarray*}
 \int_{B_{r_0}(0)\setminus B_{\frac{|x|}2}(x)} \frac{|y|^{\tau_+(\mu)} f(y)}{|x-y|^{N-2+\tau_+(\mu)}}  dy
 &\le& c\int_{B_{r_0}(0)\setminus B_{\frac{|x|}2}(x)} \frac{ |y|^{\tau_+(\mu)} f(y)}{(|x|+|y|)^{N-2+\tau_+(\mu)} }  dy
 \\&=&c|x|^{2}\int_{B_{\frac{r_0}{|x|}}(0)\setminus B_{\frac{1}2}(e_x)} \frac{ |z|^{\tau_+(\mu)} f(|x|z)}{(1+|z|)^{N-2+\tau_+(\mu)}}  dz.
\end{eqnarray*}
For any $\varepsilon>0$ fixed, there exists $R_\varepsilon>0$ such that $1/(1+R_\varepsilon)^{N-2+\tau_+(\mu)} <\varepsilon$, then
\begin{eqnarray*}
\int_{B_{\frac{r_0}{|x|}}(0)\setminus B_{R_\varepsilon}(0)} \frac{|z|^{\tau_+(\mu)} f(|x|z)}{(1+|z|)^{N-2+\tau_+(\mu)}}  dz&\le&\varepsilon
\int_{B_{\frac{r_0}{|x|}}(0)\setminus B_{R_\varepsilon}(0)}  |z|^{\tau_+(\mu)} f(|x|z)   dz
\\&\le& \varepsilon |x|^{-N-\tau_+(\mu)}\norm{f}_{L^1(\Omega,\Gamma_\mu dx)}
\end{eqnarray*}
and
\begin{eqnarray*}
\int_{B_{R_\varepsilon}(0)} \frac{|z|^{\tau_+(\mu)} f(|x|z)}{(1+|z|)^{N-2+\tau_+(\mu)}}  dz&\le&
\int_{  B_{R_\varepsilon}(0)}  |z|^{\tau_+(\mu)} f(|x|z)   dz
\\&\le&  |x|^{-N-\tau_+(\mu)}\int_{  B_{ |x| R_\varepsilon}(0)}  |y|^{\tau_+(\mu)} f(y)   dy,
\end{eqnarray*}
where from (\ref{f1}) we have that
$$\int_{  B_{ |x| R_\varepsilon}(0)}  |y|^{\tau_+(\mu)} f(y)   dy\to 0\quad  {\rm as}\quad |x|\to0.$$
So $\lim_{|x|\to0}\mathbb{E}_1(x)=0$.  Therefore, along with (\ref{4.5}) and (\ref{4.1-31}),
$$\lim_{|x|\to0}\mathbb{G}_\mu[f](x)\Phi_\mu^{-1}(x)=0.$$

 {\it Case 2:    $\mu_0<0$ and  $\mu_0\le  \mu<0$.} We observe that $N\ge 3$ and
 \begin{eqnarray*}
\mathbb{G}_\mu[f](x) &\le& c\displaystyle\left(\int_\Omega \frac{f(y)}{|x-y|^{N-2}} dy+\int_\Omega\frac{|x|^{\tau_+(\mu)} f(y)}{|x-y|^{N-2+\tau_+(\mu)}} dy\right.\\&&+\left.\int_\Omega\frac{ |y|^{\tau_+(\mu)}f(y) }{|x-y|^{N-2+\tau_+(\mu)}}  dy+ \int_\Omega\frac{|x|^{\tau_+(\mu)}|y|^{\tau_+(\mu)}}{|x-y|^{N-2+2\tau_+(\mu)}}f(y) dy\displaystyle\right).
 \end{eqnarray*}
From (\ref{4.1-31}), we have that
 \begin{eqnarray*}
 \lim_{|x|\to0}|x|^{-\tau_-(\mu)} \int_{B_{\frac{|x|}2}(x)} \frac{  f(y)}{|x-y|^{N-2 }}  dy =0.
\end{eqnarray*}
Moreover, we have that
\begin{eqnarray*}
  \int_{B_{\frac{|x|}2}(x)} \frac{|y|^{\tau_+(\mu)} f(y)}{|x-y|^{N-2+\tau_+(\mu)}}  dy &\le& \sup_{z\in B_{\frac{3|x|}2}(0)\setminus B_{\frac{|x|}2}(0) } |z|^{\tau_+(\mu)}f(z)\, \int_{B_{\frac{|x|}2}(x)} |x-y|^{2-N+\tau_+(\mu)}  dy
  \\&\le& c|x|^{2}\sup_{B_{\frac{3|x|}2}(0)\setminus B_{\frac{|x|}2}(0)} f(z),
\end{eqnarray*}
and then
\begin{eqnarray*}
 \lim_{|x|\to0}|x|^{-\tau_-(\mu)} \int_{B_{\frac{|x|}2}(x)} \frac{|y|^{\tau_+(\mu)}  f(y)}{|x-y|^{N-2 +\tau_+(\mu)}}  dy =0.
\end{eqnarray*}
 Similarly,
\begin{eqnarray*}
 \lim_{|x|\to0}|x|^{-\tau_-(\mu)} \int_{B_{\frac{|x|}2}(x)} \frac{|x|^{\tau_+(\mu)} f(y)}{|x-y|^{N-2+\tau_+(\mu)}} dy =0\quad \\{\rm and}\quad \lim_{|x|\to0}|x|^{-\tau_-(\mu)} \int_{B_{\frac{|x|}2}(x)} \frac{|x|^{\tau_+(\mu)}|y|^{\tau_+(\mu)}}{|x-y|^{N-2+2\tau_+(\mu)}}f(y)  dy =0.
\end{eqnarray*}

For $y\in B_{r_0}(0)\setminus B_{\frac{|x|}2}(x)$, we have that
$$|x-y|\ge \frac18 (|x|+|y|)$$ and
\begin{eqnarray*}
 \int_{B_{r_0}(0)\setminus B_{\frac{|x|}2}(x)} \frac{ f(y)}{|x-y|^{N-2}}  dy
 &\le& c\int_{B_{r_0}(0)\setminus B_{\frac{|x|}2}(x)} \frac{ f(y)}{(|x|+|y|)^{N-2} }  dy
 \\&=&c|x|^{2}\int_{B_{\frac{r_0}{|x|}}(0)\setminus B_{\frac{1}2}(e_x)} \frac{ f(|x|z)}{(1+|z|)^{N-2+\tau_+(\mu)}}  dz.
\end{eqnarray*}
Fixed $\varepsilon>0$, there exists $R_\varepsilon>1$ such that $1/(1+R_\varepsilon)^{N-2} <\varepsilon$, then
for $|z|>R_\varepsilon$, we have that $$\left(\frac{1+|z|}{|z|}\right)^{\tau_+(\mu)}\le c, $$
\begin{eqnarray*}
\int_{B_{\frac{r_0}{|x|}}(0)\setminus B_{R_\varepsilon}(0)}\frac{ f(|x|z)}{(1+|z|)^{N-2}}  dz&\le& c \int_{B_{\frac{r_0}{|x|}}(0)\setminus B_{R_\varepsilon}(0)} \frac{|z|^{\tau_+(\mu)} f(|x|z)}{(1+|z|)^{N-2+\tau_+(\mu)}}  dz\\&\le&\varepsilon
\int_{B_{\frac{r_0}{|x|}}(0)\setminus B_{R_\varepsilon}(0)}  |z|^{\tau_+(\mu)} f(|x|z)   dz
\\&\le& \varepsilon |x|^{-N-\tau_+(\mu)}\norm{f}_{L^1(\Omega,\Gamma_\mu dx)}
\end{eqnarray*}
and for $|z|<R_\varepsilon$, by the fact that $\frac{2-N}2\le \tau_+(\mu)<0$,
\begin{eqnarray*}
\int_{B_{R_\varepsilon}(0)} \frac{  f(|x|z)}{(1+|z|)^{N-2 }}  dz&\le& \int_{  B_{R_\varepsilon}(0)} |z|^{\tau_+(\mu)} f(|x|z)   dz
\\&\le& c |x|^{-N-\tau_+(\mu)}\int_{  B_{ |x| R_\varepsilon}(0)}  |y|^{\tau_+(\mu)} f(y)   dy,
\end{eqnarray*}
where
$$\int_{  B_{ |x| R_\varepsilon}(0)}  |y|^{\tau_+(\mu)} f(y)   dy\to 0\quad  {\rm as}\quad |x|\to0.$$
 So
$$\lim_{|x|\to0} \Phi_\mu^{-1}(x)  \int_{B_{r_0}(0)\setminus B_{\frac{|x|}2}(x)}
\frac{  f(y)}{|x-y|^{N-2 }}   dy=0,$$
then we have that
$$\lim_{|x|\to0}\left(\int_{\Omega } \frac{  f(y)}{|x-y|^{N-2 }}  dy\right) \Phi_\mu^{-1}(x)=0.$$
Similar we can prove
$$\lim_{|x|\to0}\left(\int_{\Omega} \frac{ |x|^{\tau_+(\mu)} f(y)}{|x-y|^{N-2+\tau_+(\mu) }}  dy\right) \Phi_\mu^{-1}(x)=0,$$
$$\lim_{|x|\to0}\left(\int_{\Omega} \frac{ |y|^{\tau_+(\mu)} f(y)}{|x-y|^{N-2+\tau_+(\mu) }}  dy\right) \Phi_\mu^{-1}(x)=0$$
and
$$\lim_{|x|\to0}\left(\int_{\Omega} \frac{|x|^{\tau_+(\mu)} |y|^{\tau_+(\mu)} f(y)}{|x-y|^{N-2+2\tau_+(\mu) }}  dy\right) \Phi_\mu^{-1}(x)=0.$$
Thus (\ref{1.01}) holds.\hfill$\Box$

\begin{remark}
We remark that for $\tau> \tau_-(\mu)$, if
\begin{equation}\label{a111}
  \limsup_{x\to0} f(x)|x|^{-\tau  +2}<+\infty,
\end{equation}
then $f$ verifies (\ref{f1}) and (\ref{4.1-3}). Furthermore, when $\mu>\mu_0$ and $f\in C^\gamma_{loc}(\overline{\Omega}\setminus \{0\})$ verifies (\ref{a111}),  there exists $c>0$ such   that
$$|\mathbb{G}_\mu[f](x)|\le c |x|^{\min\{\tau,\, \tau_+(\mu)\}}.$$
%if $\mu=\mu_0$, problem (\ref{eq 1.1f}) has a unique solution $u_k$ satisfying
%$$\lim_{|x|\to0}\frac{|u_k(x)-kG_\mu(x)|}{G_\mu(x)}=0.$$

\end{remark}

\subsection{Existence and nonexistence }

This subsection is devoted to build the distributional identity for the  nonhomogeneous problem (\ref{eq 1.1f}) .

\begin{lemma}\label{pr 4.1}
Assume that  $f$ is a function in $C^\gamma_{loc}(\overline{\Omega}\setminus \{0\})$ satisfying (\ref{f1}).
  Then for $\forall k\in\R$,  problem (\ref{eq 1.2}) admits a unique weak solution  $u_k$,
has a  unique  solution $u_k$, which is a classical solution of  problem (\ref{eq 1.1f}).

Furthermore, if (\ref{4.1-3}) holds true, then    $\lim_{x\to0}u_k(x)\Phi_\mu^{-1}(x)=k$.
\end{lemma}
\noindent{\bf Proof.} Let $f_n=f\eta_n$, where
$\eta_n(r)=1-\eta_0(nr)$ for $r\ge0$. We see that $f_n$ is bounded. Let  $v_n$ and $v_n^+$  be the solution of
(\ref{2.02}) replaced $f$ by $f_n$ and $|f_n|$ respectively, then
\begin{equation}\label{3.1}
 |v_n|\le v_n^+
\end{equation}
and  for any $\xi\in C^{1.1}_0(\Omega)$,
$$\int_{\Omega} v_n\,   \mathcal{L}_\mu^*(\xi )\, d\mu =\int_{\Omega} f_n  \xi \, d\mu,\qquad \int_{\Omega} v_n^+\,  \mathcal{L}_\mu^*(\xi ) d\mu=\int_{\Omega} |f_n|  \xi d\mu.$$
Take $\xi=\xi_0$ as before, we have that $ v_n^+$ is uniformly bounded in $L^1(\Omega, \,d\mu)$ by (\ref{f1}), so is $v_n$.

{\it Claim 1:  $v_n^+$ is uniform bounded in $L^1(\Omega, \,|x|^{-1} d\mu)$.}
 For  $\sigma>0$,  let $\xi_\sigma$ be defined in (\ref{t1}), then using the expression and the estimate of
 $\mathcal{L}_\mu^*(\xi_\sigma )(x)$ in $B_{\frac1{n_0}}(0)$ given by (\ref{eq 3.3000}) and (\ref{eq 3.3100}) respectively, we have that
 %\begin{equation}\label{eq 4.2000}
 %\mathcal{L}_\mu^*(\xi_\sigma )(x)=\arraycolsep=1pt\left\{
%\begin{array}{lll}
 %\frac{N-1+\tau_+(\mu)}{|x|}\quad{\rm for} \quad  \sigma\le |x|<1/n_0,\\[2mm]
 %\frac1\sigma+N-1+\tau_+(\mu)\quad{\rm for} \quad  |x|<\sigma
%\end{array}
%\right.
%\end{equation}
%and
%there exists $c>0$ such that
%$$|\mathcal{L}_\mu^*(\xi_\sigma )(x)|\le c,\qquad\forall\, x\in \Omega \setminus B_{\frac1{n_0}}(0),$$
%where $N-1+\tau_+(\mu)>0$ for $\mu\ge\mu_0$.
%Therefore, we have that
\begin{eqnarray*}
(N-1+\tau_+(\mu))\int_{B_{\frac1{n_0}}(0)\setminus B_\sigma(0)} v_n^+\,  |x|^{-1} d\mu &\le&
 \int_{B_{\frac1{n_0}}(0)} v_n^+\,  \mathcal{L}_\mu^*(\xi_\sigma ) d\mu
 \\ &\le & \int_\Omega |f| \xi_\sigma  d\mu +c\int_{\Omega\setminus B_{\frac1{n_0}}(0)} v_n^+  d\mu
 \\ &\le & c \int_\Omega |f| \rho \, d\mu
\end{eqnarray*}
Passing to the limit as $\sigma\to0$, we have that
$$\int_{B_{\frac1{n_0}}(0)} v_n^+\,  |x|^{-1} d\mu\le c'\int_\Omega |f| \rho \, d\mu,$$
which, together with the fact that
$$\int_{\Omega\setminus B_{\frac1{n_0}}(0)} v_n^+\,  |x|^{-1} d\mu\le  c\int_\Omega |f| \rho \, d\mu,$$
implies that
$$\int_{\Omega} v_n^+\,  |x|^{-1} d\mu\le c'\int_\Omega |f|\rho  \, d\mu,$$
that is,  $v_n^+$ is uniform bounded in $L^1(\Omega, \,|x|^{-1} d\mu)$.

Moreover,  $ \{v_n^+\}$ is increasing, and then there exists $v_+$ such that
$$v_n^+\to v_+\quad{\rm a.e.\ in}\ \Omega\quad{\rm and\ in} \ L^1(\Omega,\,|x|^{-1} d\mu).$$
Then we have that
$$\int_{\Omega} v_+  \mathcal{L}_\mu^*(\xi )\, d\mu =\int_{\Omega} |f|   \xi \,  d\mu,\quad\forall\, \xi\in C^{1.1}_0(\Omega).$$
Since $f\in C^\gamma(\overline{\Omega}\setminus \{0\})$, then from  Lemma 4.10 in \cite{C}, we have that
$v\in C^1(\Omega\setminus\{0\})$ and then from (\ref{3.1}), up to subsequence,  there exists $v_f$ such that
 $$v_n\to v_f \quad {\rm  \ in} \ C^1(\Omega\setminus\{0\}) \quad{\rm and\ \ in}\ \, L^1(\Omega, \,d\mu)$$
and
$$\int_{\Omega} v_f  \mathcal{L}_\mu^*(\xi ) \,d\mu =\int_{\Omega} f  \xi \,  d\mu ,\quad\forall\, \xi\in   C^{1.1}_0(\Omega),$$
thus, $u_{k,f}=kG_\mu+v_f$   is a weak solution of (\ref{eq 1.2}) and the uniqueness follows by the Kato's inequality Proposition  \ref{pr 2.1}.  \smallskip

{\it Claim 2:  $v_f$ is a classical solution of (\ref{eq 1.1f}) subject to $\lim_{x\to0}u(x)\Phi_\mu^{-1}(x)=0$.}
From Corollary 2.8 in \cite{V} with $L^*=\mathcal{L}_{\mu}^*$, which is strictly elliptic in $\Omega\setminus B_r(0)$, we have that for $q<\frac{N}{N-1}$,
\begin{eqnarray}
\norm{v_n\Gamma_\mu}_{W^{1,q}(\Omega_{2r})}  &\le & c\norm{f\Gamma_\mu}_{L^1 (\Omega\setminus B_r(0))}+ c\norm{v_n\Gamma_\mu}_{L^1(\Omega\setminus B_r(0))} \nonumber\\
  &\le & c \norm{f }_{L^1 (\Omega,\,d\mu)},\label{5.1}
\end{eqnarray}
where $\Omega_{2r}=\{x\in\Omega\setminus B_{2r}(0):\, \rho(x)>2r\}.$
 We see that
\begin{eqnarray*}
-\Delta v_n   = -\frac{\mu}{|x|^2}v_n + f.
\end{eqnarray*}
For any compact set $K$ in $\Omega$ away from the origin, it is standard to improve the regularity $v_n$
$$\norm{v_n}_{C^3(K)}\le c[\norm{f }_{L^1 (\Omega,\,d\mu)}+ \norm{f}_{C^1 (K)}]$$
where $c>0$ is independent of $n$. Then $v_f$ is a classical solution of  (\ref{eq 1.1f}) by the stability theorem, so is $kG_\mu+v_f$.

We observe that
$$v_n=\mathbb{G}_\mu[f_n]\quad{\rm  and}\quad v_f=\lim_{n\to+\infty}\mathbb{G}_\mu[f_n]=\mathbb{G}_\mu[f],$$
then it deduces from Proposition \ref{pr 4.1.1} if (\ref{4.1-3}) holds true, that
$\lim_{x\to0}v_f(x)\Phi_\mu^{-1}(x)=0$.
We complete the proof. \hfill$\Box$

 \begin{lemma}\label{lm 2.2}
Assume that   $f$ is a nonnegative function in $C^\gamma_{loc}(\overline{\Omega}\setminus \{0\})$ satisfying (\ref{f2}).
Then problem (\ref{eq 1.1f})
has no nonnegative solution.

\end{lemma}
{\bf Proof.} By contradiction, we assume that problem (\ref{eq 1.1f}) has a nonnegative solution of $u_f$.
Let $\{r_n\}_n$ be a sequence of strictly decreasing positive numbers converging to zero.
From (\ref{f2}) and the fact $f\in C^\gamma(\overline{\Omega}\setminus \{0\})$, for any $r_n$, we have that
$$
 \lim_{r\to0^+} \int_{B_{r_n}(0)\setminus B_r(0)}f(x) d\mu  =+\infty,
$$
then there exists $R_n\in (0,r_n)$ such that
$$
  \int_{B_{r_n}(0)\setminus B_{R_n}(0)}f d\mu =n,
$$
{\it The case of $\mu\ge0$.} Let $\delta_n=\frac1n \Gamma_\mu f\chi_{B_{r_n}(0)\setminus B_{R_n}(0)}$, then the problem
$$
 \arraycolsep=1pt\left\{
\begin{array}{lll}
 \displaystyle  \mathcal{L}_\mu u\cdot\Gamma_\mu= \delta_n\qquad
   {\rm in}\quad  {\Omega}\setminus \{0\},\\[1mm]
 \phantom{  L_\mu -- }
 \displaystyle  u= 0\qquad  {\rm   on}\quad \partial{\Omega},\\[1mm]
 \phantom{   }
  \displaystyle \lim_{x\to0}u(x)\Phi_\mu^{-1}(x)=0
 \end{array}\right.
$$
has a unique positive solution  $w_n$  satisfying (in the usual sense)
$$\int_{\Omega} w_n \mathcal{L}_\mu(\Gamma_\mu\xi) dx=\int_{\Omega} \delta_n \xi dx,\quad\forall\, \xi\in C^{1.1}_0(\Omega).$$
For any $\xi\in C^{1.1}_0(\Omega)$, we have that
$$\int_\Omega w_n  \mathcal{L}_\mu^*(\xi)\,  d\mu =\int_{\Omega} \delta_n \xi\, dx\to \xi(0)\quad{\rm as}\quad n\to+\infty. $$
Therefore, by Lemma \ref{lm 2.3}  for any compact set $\mathcal{K}\subset \Omega\setminus \{0\}$
$$\norm{w_n-G_\mu}_{C^1(\mathcal{K})}\to 0\quad{\rm as}\quad  {n\to+\infty}.$$
So we fixed a point $x_0\in \Omega\setminus \{0\}$, let $r_0=\frac{\min\{|x_0|,\, \rho(x_0)\}}{2}$ and $\mathcal{K}=  \overline{B_{r_0}(x_0)}$, then
there exists $n_0>0$ such that for $n\ge n_0$,
\begin{equation}\label{4.1}
 w_n\ge \frac12G_\mu\quad{\rm in}\quad \mathcal{K}.
\end{equation}

Let $u_n$ be the solution (in the usual sense) of
$$
 \arraycolsep=1pt\left\{
\begin{array}{lll}
 \displaystyle  \mathcal{L}_\mu u\cdot\Gamma_\mu= n\delta_n\qquad
   {\rm in}\quad  {\Omega}\setminus \{0\},\\[1mm]
 \phantom{  L_\mu -- }
 \displaystyle  u= 0\qquad \ \ {\rm   on}\quad \partial{\Omega},\\[1mm]
 \phantom{   }
  \displaystyle \lim_{x\to0}u(x)\Phi_\mu^{-1}(x)=0,
 \end{array}\right.
$$
then we have that
$$u_n\ge nw_n\quad {\rm in}\quad  \Omega.$$
Together with (\ref{4.1}), we derive that
$$u_n\ge  \frac n2G_\mu\quad {\rm in}\quad  \mathcal{K}.$$
Then by Comparison Principle, we have that
$$u_f(x_0)\ge u_n(x_0)\to+\infty\quad{\rm as}\quad n\to+\infty,$$
which contradicts that $u_f$ is classical solution of (\ref{eq 1.1f}). \medskip

 {\it The case of $\mu_0\le \mu<0$.}
 Let $w_n$ be the solution of
$$
 \arraycolsep=1pt\left\{
\begin{array}{lll}
 \displaystyle  \mathcal{L}_\mu u(x)=  f\chi_{B_{r_n}(0)\setminus B_{R_n}(0)}\qquad
   &\forall x\in \Omega\setminus \{0\},\\[1.5mm]
 \phantom{  L_\mu   }
 \displaystyle  u(x)= 0\qquad  &\forall x\in \partial{ \Omega},\\[1.5mm]
 \phantom{   }
  \displaystyle \lim_{x\to0}u(x)\Phi_\mu^{-1}(x)=0,
 \end{array}\right.
$$
then it follows by comparison principle
$$u\ge w_n(x)\ge \int_{B_{r_n}(0)\setminus B_{R_n}(0) }\mathcal{G}_\mu(x,y) f(y) dy.$$
It follows from  Remark  \ref{re 4.1} that for $x,y\in B_3(0)\setminus\{0\}$, $x\not=y$,
$$\mathcal{G}_\mu(x,y)\ge
c\left(|x-y|^{2-N}+\frac{|x|^{\tau_+(\mu)} }{|x-y|^{N-2+\tau_+(\mu)}}+\frac{ |y|^{\tau_+(\mu)}}{|x-y|^{N-2+\tau_+(\mu)}}+ \frac{|x|^{\tau_+(\mu)}|y|^{\tau_+(\mu)}}{|x-y|^{N-2+2\tau_+(\mu)}}\right).$$
For $x_0\in \R^N$ with $|x_0|=2$ fixed,   it deduce that
\begin{eqnarray*}
u(x_0)\ge w_n(x_0) &\ge & c2^{\tau_+(\mu)}\int_{B_2(0)} \frac{|y|^{\tau_+(\mu)}}{|x_0-y|^{N-2+2\tau_+(\mu)}}f(y) dy  \\
  &\ge &  c\int_{B_{r_n}(0)\setminus B_{R_n}(0)}f\, d\mu=cn\to+\infty\quad{\rm as}\quad n\to+\infty,
\end{eqnarray*}
which is impossible.
  \hfill$\Box$

%%%%%%%%%%%%%%%%%%%%%%%%%%%%%%%%%%%%%%%%%%%%%%%%%%%%%%%%%%%%%%%%%%%%%%%%%%%%%%%%%%%%%%%%%%%%%%%%%%%%%%%%%%%%%%%%%%%%%%%%%%%%%%%%%%%%%%%%%%%%%%%%%%%%%%%%%%%%%%%%%%%%%%%%%%%%%%%%%%%%%%%%%%%%

\subsection{Classification }
In this subsection, we are devoted to classify the isolated singular solutions of (\ref{eq 1.1f})
in the distributional sense. When $\mu=0$, the related classification of isolated singularities was studied
in \cite{BL}.

\begin{proposition}\label{lm 4.3}
Assume that  $f$ is a function in $C^\gamma_{loc}(\overline{\Omega}\setminus \{0\})$ satisfying (\ref{f1}) and
$u$ is a nonnegative solution of (\ref{eq 1.1f}).
  Then there exists some $k\ge 0$ such that  $u$ is a the $d\mu$-distributional solution of (\ref{1.2f}).
  % with such $k$.
\end{proposition}
{\bf Proof.} Let
$$ \bar u(r)=|\mathcal{S}^N|^{-1}r^{1-N}\int_{\partial B_r(0)} u(x) d\omega(x).$$
For $r\in(0,r_0)$, we have that
$$-\bar u''(r)-\frac{N-1}{r}\bar u'(r)+\frac{\mu}{r^2} \bar u(r)\ge \bar f(r),$$
where
$$\bar f(r)=r^{1-N}\int_{\partial B_r(0)} f(x) d\omega(x).$$
Denote $$\bar u(r)=r^{\tau_+(\mu)} v(r),$$
then
$$-v''(r)-\frac{N+2\tau_+(\mu)-1}{r} v'(r)\ge r^{-\tau_+(\mu)} \bar f,$$
where $N+2\tau_+(\mu)$ plays the dimensional role. By (\ref{f1}), we have that
$$\int_0^{r_0}  r^{-\tau_+(\mu)} \bar f(r) r^{N+2\tau_+(\mu)-1}dr= \int_{B_{r_0}(0)}  |x|^{-\tau_+(\mu)}|f(x)| dx  <+\infty, $$
where $r_0>0$ such that $B_{r_0}(0)\subset \Omega$.
Following the step 1 in the proof of Theorem 1.1 in \cite{BL}, there exists $c>0$ such that
$$v(r)\le  \left\{\arraycolsep=1pt
\begin{array}{lll}
c |x|^{2-(N+2\tau_+(\mu))}\quad
   &{\rm if}\quad N+2\tau_+(\mu)\ge 3,\\[1mm]
 \phantom{   }
- c\ln|x| \quad  &{\rm   if}\quad N+2\tau_+(\mu)=2,
 \end{array}
 \right.$$
that is, $\bar u\le c\Phi_\mu$.
So for $\xi\in  C^\infty_c(\Omega)$, it is well-defined that
\begin{eqnarray*}
|\int_\Omega u \mathcal{L}^*_\mu(\xi)\, d\mu| &\le &\norm{\xi}_{C^2}\int_{B_{r_0}(0)} u(x) |x|^{-1} \Gamma_\mu (x)dx +c \\
    &\le &   c\norm{\xi}_{C^2}\int_0^{r_0} \bar u(r) r^{-1+\tau_+(\mu)} r^{N-1} dr +c
    \\&\le& c\norm{\xi}_{C^2}\int_0^{r_0}r^{\tau_-(\mu)+\tau_+(\mu)+N-2}dr+c
    \\&<&+\infty.
\end{eqnarray*}

 We observe that for $\xi\in  C^\infty_c(\Omega\setminus\{0\})$, it follows by Divergence theorem that
$$\int_{\Omega}u  \mathcal{L}^*_{\mu}(\xi)   \Gamma_{\mu}dx -\int_{\Omega}f\,\xi \Gamma_{\mu} dx=0.$$
By Schwartz Theorem (\cite[Theorem XXXV]{S}),  there exists a multiple index $p$,
$$u \Gamma_{\mu} \mathcal{L}^*_{\mu}-f\Gamma_{\mu}=\sum_{|a|=0}^p k_a D^{a}\delta_0,$$
i.e. for any $\xi\in C^\infty_c(\Omega)$
\begin{equation}\label{4.2}
 \int_{\Omega}u (\mathcal{L}^*_{\mu}\xi-f\xi)  \, d\mu =\sum_{|a|=0}^{|p|} k_a D^{a}\xi(0).
\end{equation}
We are left to show that $k_a=0$ for $|a|\ge 1$.
For multiple index $\bar a\not=0$, taking $\xi_{\bar a}(x)=x_i^{\bar a_i}\eta_{n_0}$ and denoting $\xi_{\bar a,\varepsilon}(x)=\xi_{\bar a}(\frac{x}{\varepsilon})$,  we have that $\xi_{\bar a}\in C^\infty_c(\Omega)$ and then for $\varepsilon\in(0,\, \frac12)$,
\begin{eqnarray*}
  \mathcal{L}^*_{\mu} \xi_{\bar a,\varepsilon}(x)  = \frac{1}{\varepsilon^2}  (-\Delta) \xi_{\bar a}(\frac{x}{\varepsilon})-\frac{2}{\varepsilon} \frac{x  }{|x|^2}  \cdot \nabla\xi_{\bar a}(\frac{x}{\varepsilon}),
\end{eqnarray*}
and on the one side,
\begin{eqnarray*}
\left| \int_{\Omega}u  \mathcal{L}^*_{\mu}( \xi_{\bar a,\varepsilon}) d\mu \right| &=& \left| \frac{1}{\varepsilon^2} \int_{B_{2\varepsilon}(0)}u_{\mu}\Gamma_{\mu} (-\Delta)  \xi_{\bar a}(\frac{x}{\varepsilon}) dx-  \frac{1}{\varepsilon}  \int_{B_{2\varepsilon}(0)}u_{\mu}\Gamma_{\mu} \frac{x  }{|x|^2}  \cdot \nabla\xi_{\bar a}(\frac{x}{\varepsilon})\right| \\
    &\le &\frac{1}{\varepsilon^2}\int_{B_{2\varepsilon}(0)}u_{\mu}\Gamma_{\mu} dx+\frac{1}{\varepsilon}     \int_{B_{2\varepsilon}(0)} \frac{ u_{\mu}\Gamma_{\mu} }{|x|} dx
    \\&\le&   \left\{\arraycolsep=1pt
\begin{array}{lll}
c  \quad
   &{\rm if}\quad N\ge 3,\\[1mm]
 \phantom{   }
- c\ln \varepsilon \quad  &{\rm   if}\quad N=2,
 \end{array}
 \right.
\end{eqnarray*}
where $c$ is independent of $\varepsilon$. Moreover, we have that
$$\left| \int_{\Omega}u_{\mu}   \xi_{\bar a,\varepsilon}   d\mu \right|\le \norm{\xi_{\bar a}}_{L^\infty} \int_{B_{2\varepsilon}(0)} f_{\mu}\Gamma_{\mu} dx\to0\quad{\rm as} \ \ \varepsilon\to 0^+. $$
On the other side,
$$\sum_{|a|=0}^{|p|} k_a D^{a} \xi_{\bar a,\varepsilon}(0) =\frac{k_{\bar a}}{\varepsilon^{|{\bar a}|}}|{\bar a!}|,  $$
where
$$|{\bar a}|=\sum\bar a_i\quad{\rm and}\quad \bar a!=\prod_{i=1}^N(\bar a_i)!\ge1.$$
So if $k_{\bar a}\not=0$, we have that
$$\left|\sum_{|a|=0}^{|p|} k_a D^{a} \xi_{\bar a,\varepsilon}(0)\right|\to +\infty\quad{\rm as} \ \ \varepsilon\to 0^+,$$
that is, the right hand of (\ref{4.2}) with $\xi=\xi_{\bar a,\varepsilon}$ blows up with the rate $\varepsilon^{-|\bar a|}$,
while the left hand of (\ref{4.2}) keeps bounded for $N\ge 3$ or blows up in $N=2$, but controlled by $-\ln \varepsilon$ as $\varepsilon\to0^+$.
This is a contradiction and
so $k_a=0$ for $|a|\ge 1$.

Therefore, we have that
\begin{equation}\label{4.3}
 \int_{\Omega}(u_{\mu} \mathcal{L}^*_{\mu}\xi-  f\xi)\, d\mu  = k_0 \xi(0),\quad\forall\, \xi\in C^\infty_c(\Omega).
\end{equation}
For $\xi\in C^{1.1}_0(\Omega)$, take a sequence of functions in  $C^\infty_c(\Omega)$ converging to $\xi$, then the identity (\ref{4.3}) holds
for any $\xi\in C^{1.1}_0(\Omega)$.\hfill$\Box$\medskip

\noindent{\bf Proof of Theorem \ref{teo 3}.} The part $(i)$ follows by Lemma \ref{pr 4.1} and  Proposition \ref{pr 4.1.1}, the part $(ii)$ does by Proposition \ref{lm 4.3}. Lemma \ref{lm 2.2} implies the nonexistence in part $(iii)$.
 \hfill$\Box$

\setcounter{equation}{0}
\section{Generalization and Application}
\subsection{$d\mu$-distributional solution }

In this subsection, we prove existence of $d\mu$-distributional solution of (\ref{eq 1.2}) with a general nonhomogeneous term.

\begin{proposition}\label{cr 1.1}
 Assume that $\Omega$ is a smooth bounded domain containing the origin in $\R^N$, $\rho(x):={\rm dist}(x,\partial\Omega)$ and  $f$ is a measurable function satisfying
 $$
\int_{\Omega} |f|\rho\,   d\mu <+\infty.
 $$
 Then for any $k\in\R$, problem (\ref{eq 1.2}) admits a unique  $d\mu$-distributional solution.% with the singularity $\lim_{x\to0}u(x)\Phi_\mu^{-1}(x)=k$.
\end{proposition}

\noindent{\bf Proof.} By the linearity of $\mathcal{L}_\mu$, we find out the solutions $u_+$ and $u_-$ of (\ref{eq 1.2}) with the nonhomogeneous term
$f_++k_+\delta_0$ and $f_-+k_-\delta_0$, respectively, where
$a_{\pm}=\max\{\pm a,0\}$. Then $u=u_+-u_-$ is the solution of (\ref{eq 1.2}).
So we assume now that  $f$ is nonnegative and $k\ge0$. We observe that  there exists a increasing nonnegative sequence
$f_n\in C^1(\bar\Omega)$ such that
$$f_n\to f\quad {\rm as}\quad n\to\infty\quad {\rm in}\ L^1(\Omega,\, \rho d\mu).$$
From Theorem \ref{teo 3}, we have that problem (\ref{eq 1.1f}), subject to $\lim_{x\to0}u(x)\Phi_\mu^{-1}(x)=k$, has a unique solution of
$u_n$ satisfying
\begin{equation}\label{5.2}
 \int_{\Omega}u_n  \mathcal{L}_\mu^*(\xi)\, d\mu  = \int_{\Omega} f_n  \xi\, d\mu +c_\mu k\xi(0).
\end{equation}
By the Comparison Principle, we have that
$$0\le u_n\le u_{n+1}\quad \quad{\rm in}\quad \Omega\setminus \{0\}.$$
Now take $\xi=\xi_0$, the solution of (\ref{eq 2.2}), we get that
$$ \|u_n\|_{L^1(\Omega,\,d\mu)}  \le \int_{\Omega} |f|\rho\,   d\mu +ck,$$
where $c>0$ is independent of $n$. By Claim 1 in the proof of Lemma \ref{pr 4.1}, we have that
$$\int_{\Omega} u_n\,  |x|^{-1} d\mu\le c\int_\Omega |f|\rho \,d\mu.$$
Since $u_n$ is increasing and uniformly bounded in $L^1(\Omega,\,  |x|^{-1}d\mu)$,
by Lebesgue's monotone convergence theorem, we have that
$u(x):=\lim_{n\to+\infty} u_n$ is in $L^1(\Omega,\,|x|^{-1}d\mu)$ and satisfies
$$
 \int_{\Omega}u  \mathcal{L}_\mu^*(\xi)\, d\mu  = \int_{\Omega} f_n  \xi\, d\mu +c_\mu k\xi(0),\qquad \forall\, \xi\in C^{1.1}_0(\Omega).
$$

The uniqueness follows by Kato's inequality directly. Indeed, let $u_1,\, u_2$ be two solutions of (\ref{eq 1.2}), then $w:=u_1-u_2$ is a $d\mu$-distributional solution of (\ref{homo}) with
$f=0$,  taking $\xi=\xi_0$ in (\ref{sign}), $\xi_0$ being the solution of (\ref{eq 2.2}), it follows that $\int_\Omega|w|dx=0$, then $u_1=u_2$.
 This ends the proof.\hfill$\Box$

\subsection{Nonexistence on $\mu<\mu_0$}

 Although we have supposed  $\mu \ge \mu_0$ throughout this paper at the beginning, we can still deal the case of $\mu<\mu_0$ and
 %In this subsection, we
 obtain an Liouville Theorem in this case.

\begin{proposition}\label{sub}

Assume that $\mu<\mu_0$ and $f$ is a measurable nonnegative function,
then problem (\ref{eq 1.1f}) has no nontrivial nonnegative solutions.
\end{proposition}
{\bf Proof.} By contradiction, we assume that $u_0$ is a nontrivial nonnegative solution of (\ref{eq 1.1f}). Then there exist $x_0\in\Omega\setminus\{0\}$, $r\in(0,|x_0|)$ and $\epsilon_0>0$ such that
$B_{2r_0}(x_0)\subset \Omega\setminus\{0\}$ and
$$u_0\ge \epsilon_0 \quad{\rm in}\quad B_r(x_0).$$
We observe that
\begin{equation}\label{5.-1}
 \mathcal{L}_{\mu_0} u_0= (\mu_0-\mu)\frac{u_0}{|x|^2}+f\ge (\mu_0-\mu)\epsilon_0 \frac{\chi_{B_{r_0}(x_0)}}{|x|^2},
\end{equation}
where $\mu_0-\mu>0$,  $\chi_{B_{r_0}(x_0)}=1$ in $B_{r_0}(x_0)$ and $\chi_{B_{r_0}(x_0)}=0$ otherwise.

When $N\ge 3$,  from Remark \ref{re 4.1} with $K=\overline{B_{r_0}(x_0)\cup B_{r_0}(0)}$, for $x\in B_{r_0}(0)\setminus \{0\}$,
\begin{eqnarray*}
u_0(x) &\ge& (\mu_0-\mu)\epsilon_0 \mathbb{G}_{\mu_0}[\chi_{B_{r_0}(x_0)}] \\
    &\ge&   (\mu_0-\mu)\epsilon_0 c|x|^{-\frac{N-2}{2}}\int_{B_{r_0}(x_0)} \frac{1}{|x-y|^{\frac{N-2}{2}}}dy
    \\&\ge& c_0|x|^{-\frac{N-2}{2}},
\end{eqnarray*}
 where $c_0=(\mu_0-\mu)\epsilon_0 c\min_{x\in\Omega}\int_{B_{r_0}(x_0)} \frac{1}{|x-y|^{\frac{N-2}{2}}}dy>0.$
Go back to (\ref{5.-1}), we see that
\begin{eqnarray*}
\int_{\Omega\setminus B_r(0)}[(\mu_0-\mu)\frac{u_0}{|x|^2}+f] d\mu &\ge& c_0\int_{B_{r_0}(0)\setminus B_r(0)}|x|^{-N}dx  \\
   &\to& +\infty\quad{\rm as}\quad r\to0^+.
\end{eqnarray*}
Apply Theorem \ref{teo 3}$(iii)$, we obtain that
\begin{equation}\label{5.-2}
\arraycolsep=1pt\left\{
\begin{array}{lll}
 \displaystyle   \mathcal{L}_\mu u=(\mu_0-\mu)\frac{u_0}{|x|^2}+f\qquad
   {\rm in}\quad  {\Omega}\setminus \{0\},\\[1.5mm]
 \phantom{   L_\mu   }
 \displaystyle  u= 0\qquad  {\rm   on}\quad \partial{\Omega}
 \end{array}\right.
\end{equation}
 has no nonnegative solution, which contradicts that $u_0$ is a nonnegative solution (\ref{5.-2}).
 So we conclude that  problem (\ref{eq 1.1f}) has no nontrivial nonnegative solutions.\smallskip

When $N=2$ and $\mu_0=0$£¬ we have that for any compact set $K$ in $\Omega$, we have that $\mathcal{G}_0(x,y)\ge c'$ for $x,y\in K$, $x\not=y$,
\begin{eqnarray*}
u_0(x) &\ge& (\mu_0-\mu)\epsilon_0 \mathbb{G}_{\mu_0}[\chi_{B_{r_0}(x_0)}] \\
    &\ge&   (\mu_0-\mu)\epsilon_0 c \int_{B_{r_0}(x_0)} (-\ln |x-y|)dy
    \\&\ge& (\mu_0-\mu)\epsilon_0 c  |B_{r_0}(x_0)|c'.
\end{eqnarray*}
Then a contradiction could be obtained   as the case of $N\ge3$.\hfill$\Box$

 \subsection{Nonexistence of the principle eigenvalue }

 In this subsection, our aim is to study the nonexistence of principle eigenvalue of the problem
 \begin{equation}\label{6.1}
\arraycolsep=1pt\left\{
\begin{array}{lll}
 \displaystyle  -\Delta u=\lambda Vu\qquad
   {\rm in}\quad  {\mathcal{O}},\\[1.5mm]
 \phantom{    }
 \displaystyle  u\in  H^1_0(\mathcal{O}),
 \end{array}\right.
\end{equation}
where  $\mathcal{O}$ is an smooth open set containing $\{0\}$ in $\R^N$ and $V$ is a measurable function.

It was proved in \cite{SW} (Theorem 2.5 ) that problem (\ref{6.1}) admits a simple principle eigenvalue under the following
hypotheses $(V_1)$ and $(V_2)$:
\begin{itemize}
\item[]
\begin{enumerate}\item[$(V_1)$]
$V\in L^p_{loc}(\mathcal{O}\setminus\{0\})$ with $p>\frac N2$;
\end{enumerate}
\begin{enumerate}\item[$(V_2)$]
Let
\begin{equation}\label{6.2}
 V^+:=\max\{0,V\}\not=0,\quad V^+=V_1+V_2,
\end{equation}
where $ V_1\in L^{\frac N2}(\mathcal{O})$  and for any $y\in\bar \mathcal{O},$
$$\lim_{x\in\mathcal{O},\, x\to y} |x-y|^2V_2(x)=0,\quad \lim_{x\in\mathcal{O},\, x\to \infty} |x|^2V_2(x)=0\
({\rm if}\ \mathcal{O}\ {\makebox{ is unbounded}}).$$
\end{enumerate}
\end{itemize}
 More references on the principle eigenvalue with indefinite potential also could be seen in  \cite{BCF,Ch1,Sw,WW}.

Our aim here is to the nonexistence of the principle eigenvalue if  $V$ is  the Hardy-Leray potential in bounded domain. We will see
that the assumptions on $V$ above seem to be "optimal" in some sense.
More precisely, our result can be stated as follows:

\begin{theorem}\label{to 5.1}
Assume that $\mathcal{O}=\Omega$ is a bounded $C^2$ domain containing the origin,  the potential $V\in C^\gamma_{loc}(\bar\Omega\setminus\{0\})$  verifies
$$ V(x)\ge \frac{a_0}{|x|^2},\quad\forall\, x\in\Omega\setminus\{0\} $$
for some constants  $\gamma\in(0,1)$ and $a_0>0$.
Then problem (\ref{6.1}) has no  principle eigenvalue.
\end{theorem}
{\bf Proof. }  By contradiction, we assume that   problem (\ref{6.1}) has  a  principle eigenvalue $\lambda_1$, and denote by $\varphi_1$
the nonnegative corresponding eigenfunction.

Since $V$ is nonnegative and $V\in C^\gamma(\bar\Omega\setminus\{0\})$, it is standard to prove that $\lambda_1>0$ and $\varphi_1\in  H^1_0(\Omega)$ is a nonnegative classical solution of
 \begin{equation}\label{6.3}
\arraycolsep=1pt\left\{
\begin{array}{lll}
 \displaystyle  -\Delta u=\lambda_1 V u\qquad
   &{\rm in}\quad  {\Omega}\setminus\{0\},\\[1.5mm]
 \phantom{  -\Delta  }
 \displaystyle  u=0\qquad &{\rm on}\quad \partial\Omega.
 \end{array}\right.
\end{equation}
So $\varphi_1$ a super solution of
 \begin{equation}\label{6.4}
\arraycolsep=1pt\left\{
\begin{array}{lll}
 \displaystyle  -\Delta u=\lambda_1a_0  \frac{ u}{|x|^2}\qquad
   &{\rm in}\quad  {\Omega}\setminus\{0\},\\[1.5mm]
 \phantom{  -\Delta  }
 \displaystyle  u=0\qquad &{\rm on}\quad \partial\Omega.
 \end{array}\right.
\end{equation}
Let $\psi_1$ be the first positive eigenfunction of $-\Delta u=\lambda u$ in $H^1_0(\Omega\setminus B_r(0))$ with $r>0$ small.
Then there exists $\epsilon\in(0,1)$ such that $\epsilon\psi_1$ is a sub solution of (\ref{6.4}). By the super and sub solutions methods,
problem (\ref{6.4}) has a nonnegative nontrivial solution $v_1\in H^1_0(\Omega)$.

In the case that $\lambda_1a_0\in (0,-\mu_0]$, it follows by Theorem \ref{teo 3} $(iii)$ with $f=0$, we have that
$v_1=kG_{-\lambda_1a_0}$ for some $k\ge0$. If $k=0$, $v_0=0$ is trivial, which is impossible.
If $k>0$, we know that $G_{-\lambda_1a_0}\not\in H^1_0(\Omega)$, which contradicts $v_1\in H^1_0(\Omega)$.

In the case that $\lambda_1a_0>-\mu_0$, it follows by Proposition \ref{sub} that there is no nonnegative nontrivial solution.
So we obtain a contradictions. This ends the proof. \hfill$\Box$

\medskip

\bigskip

 \noindent{\bf Acknowledgements:}  The authors  would like to thank Prof. Brezis for useful discussing.   
 H. Chen  is supported by NNSF of China, No:11726614, 11661045,   by SRF for ROCS, SEM and by the Jiangxi Provincial Natural Science Foundation, No: 20161ACB20007.
  A. Quaas is partially supported by  Fondecyt Grant No. 1151180 Programa Basal, CMM, U. de Chile and Millennium Nucleus Center for Analysis of PDE NC130017. F. Zhou is partially supported by NNSF of China, No:11271133, 11431005
 and Shanghai Key Laboratory of PMMP.

\end{document}